\documentclass[preprint]{elsarticle}
\usepackage[utf8]{inputenc}

\usepackage{float}
\usepackage[margin=1in]{geometry}
\usepackage{comment}
\usepackage{stackengine}
\usepackage{lmodern}
\usepackage[normalem]{ulem}

\usepackage[thinlines]{easytable}
\usepackage[ruled]{algorithm2e}
\usepackage{algpseudocode}
\usepackage{svg}
\usepackage{subcaption}
\usepackage{graphicx}
\usepackage{epsfig}
\usepackage{epstopdf}

\usepackage{amssymb}

\usepackage{bm}
\usepackage{amsmath}
\usepackage{amsfonts}
\usepackage{physics}
\usepackage{xfrac}
\usepackage[bb = dsserif]{mathalpha}
\usepackage[mathscr]{euscript}

\usepackage[]{bm}

\usepackage[breaklinks, colorlinks=false,linkcolor=blue, citecolor=blue, urlcolor=blue]{hyperref}
\usepackage[capitalise]{cleveref} 
\usepackage{doi} 
\biboptions{sort}
\usepackage[page,toc]{appendix} 
\crefname{appsec}{Appendix}{Appendix}

\begin{document}
\begin{frontmatter}

\title{Spacetime Wavelet Method for Linear Boundary-Value Problems in Sylvester Matrix Equation Form}
\author{Cody D. Cochran}
\author{Karel Matou\v{s}\corref{cor1}}
\cortext[cor1]{Corresponding author}
\journal{Journal}

\address{Department of Aerospace and Mechanical Engineering, University of Notre Dame, Notre Dame, IN 46556, USA}

\begin{abstract}
\noindent
We present a high-order spacetime numerical method for discretizing and solving linear initial-boundary value problems using wavelet-based techniques with user-prescribed error estimates. The spacetime wavelet discretization yields a system of algebraic equations resulting in a Sylvester matrix equation. We solve this system with a Global Generalized Minimal Residual (GMRES) method in conjunction with a wavelet-based recursive algorithm to improve convergence. We perform rigorous verification studies using linear partial differential equations (PDEs) with both convective and diffusive terms. The results of these simulations show the high-order convergence rates for the solution and derivative approximations predicted by wavelet theory. We demonstrate the utility of solving the Sylvester equation through comparisons to the commonly-used Kronecker product formulation. We show that our recursive wavelet-based algorithm that generates initial guesses for the iterative Global GMRES method improves the performance of the solver.
\end{abstract}


\end{frontmatter}

\section{Introduction}
\label{sec:intro}


Engineers, scientists and mathematicians have studied partial differential equations (PDEs) and techniques to find their solutions for centuries \cite{green1889essay, hadamard1908theorie}. Many of these equations do not have closed-form solutions. In the field of computational science and engineering, it is crucial to develop techniques that solve equations accurately and efficiently. Many important advancements in the field of numerical methods have been made to accomplish this goal \cite{courant1928partiellen, crank1947practical, godunov1959finite}. A subset of these long-studied equations are known as linear PDEs. Linear PDEs appear in many contexts including but not limited to fluid mechanics \cite{crank1947practical}, solid mechanics \cite{chen2006meshless}, and control theory \cite{balas2003trends}. Many different methods have been developed to discretize and solve these systems, offering certain benefits over others. In this work, we present a high-order spacetime wavelet method for solving linear PDEs. 

Multiresolution Analysis (MRA) is the mathematical foundation for wavelet theory, which builds a nested series of approximations at varying resolution levels and provides \textit{a priori} error estimates \cite{harnish2018adaptive, vasilyev2000second, dubos2013conservative}.  Wavelets are often used as tools for signal processing and image compression \cite{daubechies1992ten}. More recently, research shows that wavelets are effective tools that can be used to solve PDEs \cite{bacry1992wavelet, goedecker2009wavelets, beylkin1997adaptive, vasilyev2000second, bertoluzza1996wavelet, harnish2018adaptive, harnish2021multiresolution, harnish2023adaptive, cochran2025spacetime}. 

Despite the use of wavelets for solving PDEs, few approaches leverage wavelet-based derivatives, and instead opt for alternative derivative approximations, typically finite differences \cite{paolucci2014wamr1, paolucci2014wamr2, alam2006simultaneous, vasilyev2000second}. Solvers with non-wavelet-based derivative approximations have been proven to provide accurate solutions. However, these techniques forfeit some of the advantages provided by wavelet theory and often require a large number of nonessential grid points. Wavelet theory provides mathematically rigorous error estimates for both the solution and the derivative, both of which are known to the user \textit{a priori} \cite{goedecker2009wavelets, harnish2023adaptive, jameson1993wavelet, vasilyev2000second}. Moreover, it has been shown that wavelet derivatives demonstrate superconvergence properties \cite{jameson1993wavelet, gomes1995convergence} and have a natural extension to differentiation on sparse grids \cite{harnish2018adaptive, harnish2023adaptive}.

Most PDE solution techniques require a semi-discretization step, which produces a system of ordinary differential equations (ODEs) \cite{tsai2014two,moin2010fundamentals}. From here a decision about how to integrate in time must be carefully made to avoid computational inefficiency and unwanted sources of error. A typical choice is either an explicit or implicit time-stepping method, each of which can become inefficient under certain circumstances. The size of the time step required for the stability of explicit methods can be very small and sacrifice computational efficiency \cite{arabatzis2007implicit, merkle1988computation}. Implicit timestepping schemes can be unconditionally stable for linear problems and permit larger timesteps; however, they require a large number of computations per timestep \cite{loffeld2013comparative}. Parallel time integration techniques have been developed to improve performance by asynchronously solving time subdomains \cite{gander2007analysis, gander201550, subber2016asynchronous, prakash2014computationally, benevs2010asynchronous, subber2017asynchronous}. Although traditional time integration techniques yield solutions to systems of ODEs with local truncation error estimates expressed as a function of the time step, $\mathcal{O}(\Delta t^n)$. However, other sources of error can cause the global solution error to deviate \cite{verwer1996explicit, loffeld2013comparative}. 

An alternative to time-stepping methods is what is referred to in the literature as a spacetime method \cite{alam2006simultaneous, hughes1996space, schwab2009space}. Spacetime methods discretize and solve the PDE system in both spatial and temporal dimensions simultaneously. Much of the published spacetime work utilizes the finite element method (FEM) to discretize the problem \cite{hulbert1990space, abedi2010adaptive, takizawa2012space, hughes1996space, petersen2009space, naudet2024space}. Others use wavelet-based methods \cite{alam2006simultaneous, gunzburger2011space, schwab2009space, cochran2025spacetime}.

The traditional numerical discretization of linear elliptic or implicit time-dependent PDEs leads to systems of algebraic equations, expressed in the form $K\vec{x}=\vec{r}$, where $K$ is the Jacobian matrix, $\vec{x}$ is the solution vector, and $\vec{r}$ is the right-hand side vector. The characteristics of this discretized system gives information about which numerical methods are best suited to solve it. Direct methods such as Gaussian elimination and LU factorization are robust and can solve systems with a predictable number of steps given that $K$ is non-singular. However when dealing with large and/or sparse systems, direct methods can be inefficient \cite{saad2019iterative, glaser2009new}. Iterative methods such as the Generalized Minimal Residual (GMRES) and Conjugate Gradient (CG) methods solve the system by reducing the error at each iteration, moving the approximation towards the solution. These methods are better-equipped to efficiently handle large, sparse systems, but their convergence is dependent on a good initial guess \cite{tromeur2006choice, gao2015sensitivity, ye2020improving, saad2019iterative}. One must also be mindful to select a solver that is compatible with the eigenvalues and structure of the coefficient matrices, \textit{e.g.}, the CG method is only applicable when $K$ is symmetric and positive definite \cite{nazareth2009conjugate}.

Linear problems discretized with a spacetime wavelet method result in the Sylvester matrix equation, $AX+XB=C$, with a solution matrix $X$ and nonsymmetric $A$ and $B$ coefficient matrices. Equations of this form appear in the fields of control theory \cite{gardiner1992solution}, machine learning \cite{chen2008semi}, and PDE solution methods \cite{harker2015sylvester}, among others. Sylvester matrix equations are known to have unique solutions when the $A$ and $B$ matrices are square and the spectra of $A$ and $-B$ are disjoint, \textit{i.e.}, $A$ and $B$ do not share a common eigenvalue \cite{datta1991arnoldi}. A method for solving the Sylvester equation that often appears in the literature involves the Kronecker product and matrix vectorization to convert the system into a typical $K\vec{x}=\vec{r}$ form \cite{chiang2012sylvester, du2018fasten, chen2012projection}. Formulating and solving the Sylvester equation in this manner can become computationally inefficient as problem size increases \cite{van2000ubiquitous}. Variations of GMRES \cite{saad1986gmres} and other Krylov subspace methods are effective options when solving such systems without using the Kronecker product formulation \cite{bouhamidi2008note, jbilou1999global, robbe2002convergence, el2002block}. One popular choice of Krylov subspace method commonly used to solve the Sylvester equation is the Block-GMRES method (BGMRES). While BGMRES has an upper limit on the number of iterations required for convergence, it becomes expensive as problem size increases \cite{jbilou1999global}. 

For examples in this work, we elect to use the Global GMRES (Gl-GMRES) method developed in 1999 by Jbilou et al. The Gl-GMRES method is capable of solving the Sylvester equation and handling large, sparse, nonsymmetric $A$ and $B$ matrices with complex eigenvalues, while improving performance and decreasing memory requirements via the Modified Global Arnoldi process with an embedded modified Gram-Schmidt technique \cite{jbilou1999global, sadkane1993block}. This technique allows for a reduction in the number of orthogonal vectors/matrices generated during the solution process, reducing memory requirements and computational effort by defining a restart parameter $m$. The restarted matrix version of the Arnoldi process is referred to as the Modified Global Arnoldi process \cite{jbilou1999global, sadkane1993block}. The primary distinction between the standard Arnoldi and Global Arnoldi iteration procedures is the generation of orthogonal matrices instead of orthogonal vectors. Like other iterative methods, the convergence of the Gl-GMRES method can be improved by intelligently selecting the initial guess \cite{elleithy2001iterative, tromeur2006choice} and/or by using optimal preconditioning \cite{bouhamidi2011convex, bouhamidi2008note}. For this reason, we have developed a recursive algorithm that utilizes wavelet synthesis to generate initial guesses that improve solver performance.

This work presents a novel approach to formulate and solve linear PDEs with a spacetime wavelet solver. Example problems are discretized using entirely wavelet-based approximations, resulting in Sylvester matrix equations. Boundary and initial conditions are enforced using permutation matrices such that the system is well-posed \cite{harker2015sylvester}. The algebraic system of equations is then solved using the Gl-GMRES method with an embedded Modified Global Arnoldi algorithm. This technique avoids the computationally expensive transition to the $K\vec{x}=\vec{r}$ form. We demonstrate the capability of the solver to achieve high-order convergence rates with \textit{a priori} error estimates for both solution and derivative approximations. We show that the Gl-GMRES Sylvester formulation outperforms a restarted GMRES method when solving the Kronecker product system. We also develop a wavelet-based recursive algorithm that improves solver performance by generating informed initial guesses. The remainder of this paper is as follows: Section \ref{sec:wavelet} describes relevant wavelet theory, Section \ref{sec:comp} details the implementation of the discretization and solution techniques, and Section \ref{sec:numerical} displays and analyzes the numerical results of the spacetime wavelet solver.

\section{Wavelet Theory}
\label{sec:wavelet}
 In this section, we briefly summarize the fundamentals of spacetime wavelet theory. More detailed information on wavelets can be found in \cite{goedecker2009wavelets, donoho1992interpolating, beylkin1997adaptive, harnish2023adaptive}.
\subsection{Multiresolution Analysis}
A multiresolution analysis lays the mathematical foundation for wavelet theory. It describes sequential approximation spaces $\mathbf{V}_j$ and dual spaces $\tilde{\mathbf{V}}_j$ along with their complementary wavelet spaces and duals, denoted by $\mathbf{W}_j$ and $\tilde{\mathbf{W}}_j$, respectively. These spaces are nested within their corresponding spaces at higher resolution, having 
 properties 
\begin{align}
    \mathbf{V}_j & \subset \mathbf{V}_{j+1}, & \tilde{\mathbf{V}}_j &\subset \tilde{\mathbf{V}}_{j+1}, & \mathbf{V}_{\infty} &= L^2(\Omega), & \mathbf{V}_{j+1} = \mathbf{V}_{j} \oplus \mathbf{W}_j,
\end{align}
where $j \in \mathbb{Z}$ is the resolution level and $\Omega$ denotes some finite domain. These nested spaces provide the mathematical backbone for wavelet families of basis functions. The scaling functions ${}^0\phi^j_k$  and wavelet functions ${}^{\lambda}\psi^j_k$ are the basis functions in the $\mathbf{V}_j$ and $\mathbf{W}_j$ spaces, respectively, while their dual functions ${}^0\tilde{\phi}^j_k$ and ${}^{\lambda}\tilde{\psi}^j_k$ are the basis functions in $\tilde{\mathbf{V}}_j$ and $\tilde{\mathbf{W}}_j$. The index $k\in \mathbb{Z}$ is a local index that partially defines a location in $\Omega$ and $\lambda \in \mathbb{Z}$ describes the type and location of a wavelet function. These basis functions are represented in multiple dimensions as a tensor product of the one-dimensional components,
\begin{align}
    {}^0\Psi^j_{\vec{k}}(x_1,x_2) &= {}^0\phi^j_{k_1}(x_1) \otimes {}^0\phi^j_{k_2}(x_2).
\end{align}
In this work, we use the so-called Deslauriers-Dubuc wavelet family for the interior of the domain \cite{de2003dubuc, goedecker2009wavelets, harnish2018adaptive, harnish2021multiresolution} and boundary-modified wavelets along the domain boundary \cite{vasilyev2000second}.

\subsection{Spacetime Wavelet Discretization}
In this section, we describe the spacetime wavelet discretization as originally introduced for nonlinear problems in \cite{cochran2025spacetime}. We define a finite spacetime domain $\Omega  = \Omega_x \ \times \ [0, T] $ where $\Omega_x \subset \mathbb{R}^N$ with spatial boundary $\partial \Omega_x \subset \mathbb{R}^N$ and $t \in [0, T]$, where $[0, T] \subset \mathbb{R}^{+}$. The spacetime examples presented in this work are 2D, with one spatial dimension ($N=1$). For brevity, we consolidate all locations $k_x$, $k_t$ in the domain into an array $\vec{k}$. The spatial and temporal basis orders, $p_x$ and $p_t$, respectively, govern the size and values of this array. Similarly, we define an array that contains the spatial and temporal values: $\vec{q} = [x, t]$. The 2D spacetime wavelet representation of a function is expressed as
\begin{align}
    f^j(x, t) = f^j\big(\vec{q}\big) = \sum_{k_i \in [0, 2p_i]} {}^0\mathbb{d}^1_{\vec{k}} \ {}^0\Psi^1_{\vec{k}}(\vec{q}) + \sum_{j=1}^{j_{\mathrm{max}}} \sum_{\lambda = 1}^{2^{N+1}-1} \sum_{\vec{k}}{}^{\lambda}\mathbb{d}^j_{\vec{k}} \ {}^{\lambda}\Psi^j_{\vec{k}}(\vec{q}),
    \label{eq:2dFullapprox}
\end{align}
where $j_{\mathrm{max}}$ is the highest resolution level. 

Functions are expressed on a sparse wavelet grid by thresholding the ${}^{\lambda}\mathbb{d}$ wavelet coefficients based on some user-defined tolerance \cite{harnish2023adaptive}. However, for the examples presented in this work, we use an equivalent dense scaling function representation \cite{goedecker2009wavelets} for function approximation:
\begin{align}
    f^j\big(\vec{q}\big) = \sum_{k_i \in [0, 2^jp_i]} {}^0\mathbb{d}^j_{\vec{k}} \ {}^0\Psi^j_{\vec{k}}\big(\vec{q}\big).
    \label{eq:2dapprox}
\end{align}
This expression requires the computation of only the non-thresholded ${}^0\mathbb{d}$ coefficients to approximate the function.
Taking advantage of the properties of the Deslauriers-Dubuc wavelet family \cite{de2003dubuc}, the ${}^0\mathbb{d}$ coefficients are computed exactly by evaluating the integral with the dual basis
\begin{align}
    {}^0\mathbb{d}^j_{\vec{k}} &= \int_{\Omega} f\big(\vec{q}\big)\tilde{\Psi}^j_{\vec{k}}\big(\vec{q}\big)\mathrm{d}\Omega.
    \label{eq:dCoef}
\end{align} 

The forward wavelet transform (FWT) maps the field value coefficients ($f(\vec{q})$) to their corresponding coefficients in the wavelet space (${}^{\lambda}\mathbb{d}$) and is referred to as wavelet analysis. The inverse process of transforming wavelet coefficients to their field values is known as the backward wavelet transform (BWT), or wavelet synthesis. The operators used to perform the forward and backward transforms are denoted by $\mathbbb{F}$ and $\mathbbb{B}$ where $\mathbbb{B} = \mathbbb{F}^{-1}$. We express the discrete forward and backward wavelet transforms in index notation
\begin{align}
    {}^\lambda\mathbb{d}_{ns} &= \mathbb{F}_{no} \mathcal{F}_{or} \mathbb{F}_{rs}, & \mathcal{F}_{ns} &= \mathbf{\mathbb{B}}_{no} {}^\lambda\mathbb{d}_{or} \mathbb{B}_{r s},
    \label{eq:analSyn}
\end{align}
where $\mathcal{F}_{ns} $ is the array of function values in physical space, populated by evaluating $f(\vec{q})$ at all points on the spacetime grid. The index $n$ describes the spatial location $x$ and $s$ describes the time $t$. The field value coefficients ${}^0\mathbb{d}$ are equal to the values of $f$ evaluated at each point on the spacetime grid \cite{harnish2018adaptive, harnish2021multiresolution, goedecker2009wavelets}. Therefore, we will remove the subscript $0$ for the remainder of this work for succinctness.

\subsection{Wavelet Derivative in 2D}
In this section, we describe wavelet-based derivatives. We proceed by describing a general $\alpha$th-order derivative of a function $f(\vec{q})$ taken with respect to a spacetime variable $q_i$, \textit{i.e.}, $\frac{\partial^{(\alpha)} f(\vec{q})}{\partial q_i^{(\alpha)}}$. Similar to the approximation of a function, we express the derivative as a weighted sum
\begin{align}
     \frac{\partial^{(\alpha)}f\big(\vec{q}\big)}{\partial q_i^{(\alpha)}}\approx \sum_{\vec{k}} \mathbb{d}^j_{\vec{k}} \frac{\partial^{(\alpha)}\Psi^j_{\vec{k}}\big(\vec{q}\big)}{\partial q_i^{(\alpha)}} = \sum_{\vec{l}} \Gamma^{\vec{l}, j}_{\vec{k}}\Psi_{\vec{l}}^j\big(\vec{q}\big).
\end{align}
\noindent The Deslauriers-Dubuc wavelet basis function $\Psi\big(\vec{q}\big)$ is differentiable, therefore, we express its derivative with another wavelet approximation with so-called connection coefficients $\Gamma$ \cite{beylkin1997adaptive}. As in Eq. \ref{eq:dCoef}, we integrate with the dual basis function
\begin{align}
    \Gamma^{\vec{l}, j}_{\vec{k}} & =\int\left[\frac{\partial}{\partial q_i}\Psi_{\vec{k}}^j\big(\vec{q}\big)\right]\tilde{\Psi}_{\vec{l}}^j\big(\vec{q}\big)\mathrm{d}q_i,
    \label{eq:conn}
\end{align} 
where $l \in \mathbb{Z}$ partially describes the location of the connection coefficients on the grid.
Eq. (\ref{eq:conn}) reduces to an eigenvector problem \cite{goedecker2009wavelets} that is solved to find the values of the connection coefficients. Because the basis functions have limited differentiability, it is crucial that the interpolation orders of the bases $p_x$ and $p_t$ are large enough to accommodate the derivative orders present in the PDE. A study on the continuity and differentiability of the Deslauriers-Dubuc wavelet family is presented in \cite{rioul1992simple}.

We construct discrete operators that contain all relevant connection coefficients, denoted by ${}^{( \alpha, q_i)}\boldsymbol{\Gamma}$. These operators are defined by the domain $q_i$, resolution level $j$, interpolation order $p_i$, and derivative order $\alpha$. As we reference many times in this work, a dense wavelet derivative approximation has error described by
\begin{align}
    \bigg| \bigg| \frac{\partial^{\alpha} f\big(\vec{q}\big)}{\partial q_i^{\alpha}} - \frac{\partial^{\alpha} f^j\big(\vec{q}\big)}{\partial q_i^{\alpha}}\bigg| \bigg|_{\infty} \leq \mathcal{O}\left(\Delta q_i^{p_i-\alpha}\right),
    \label{eq:derErr}
\end{align}
where $\Delta q_i$ is the grid spacing in the $q_i$ direction \cite{dubos2013conservative, mccormick1994wavelet, harnish2018adaptive}. When solving PDEs, the expected order of convergence is defined by the smallest $p-\alpha$ value present in the PDE. For example, if we have a second-order PDE with a first-order temporal derivative ($\alpha_t = 1$) and a second-order spatial derivative ($\alpha_x = 2$)  with $p_t=4$ and $p_x=6$, the convergence rate will be dictated by the temporal derivative as $(p_t-\alpha_t)=3 < (p_x - \alpha_x) = 4$. Therefore, $3$rd-order convergence is predicted for this example. We note that there is a similar expression for function approximations using wavelets that can be found in \cite{harnish2018adaptive, harnish2023adaptive}.

\section{Computational Implementation}
\label{sec:comp}
We implement our spacetime solver computationally using the Multiresolution Wavelet Toolkit (MRWT) \cite{harnish2018adaptive, harnish2021multiresolution, harnish2023adaptive, cochran2025spacetime}. This toolkit has been developed to discretize and solve systems with wavelets using a variety of techniques. MRWT utilizes both MPI and OpenMP to create a hybrid parallelization to efficiently construct and implement wavelet operators. For matrix storage and operations, we use the Eigen C++ template library \cite{eigenweb}.

\subsection{Linear Diffusion}\label{sec:ld}
\begin{figure}[H]
\begin{subfigure}[t]{50mm}
  \includegraphics[width=50mm]{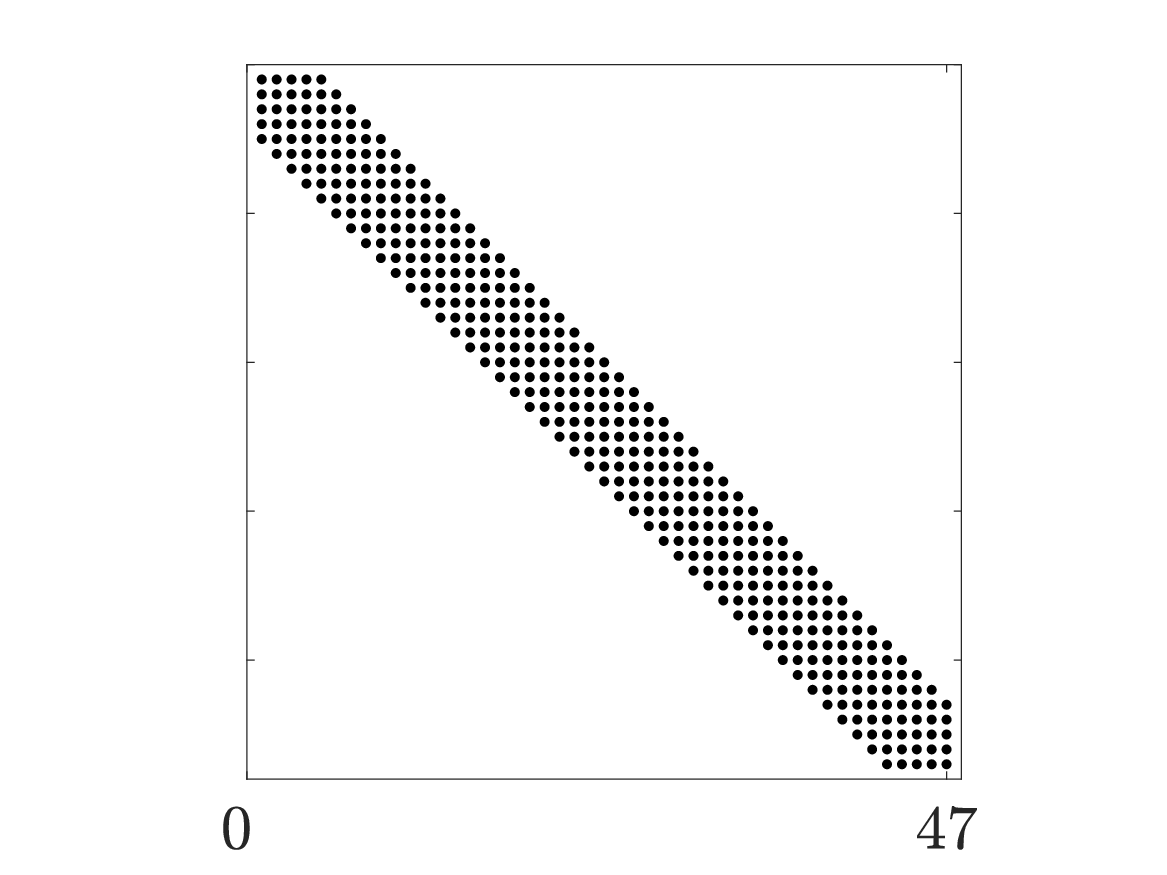}
  \caption{Sparsity pattern of the $A$ matrix.}
  \label{fig:aSparsity}
\end{subfigure}
\hfill
\begin{subfigure}[t]{50mm}
  \includegraphics[width=50mm]{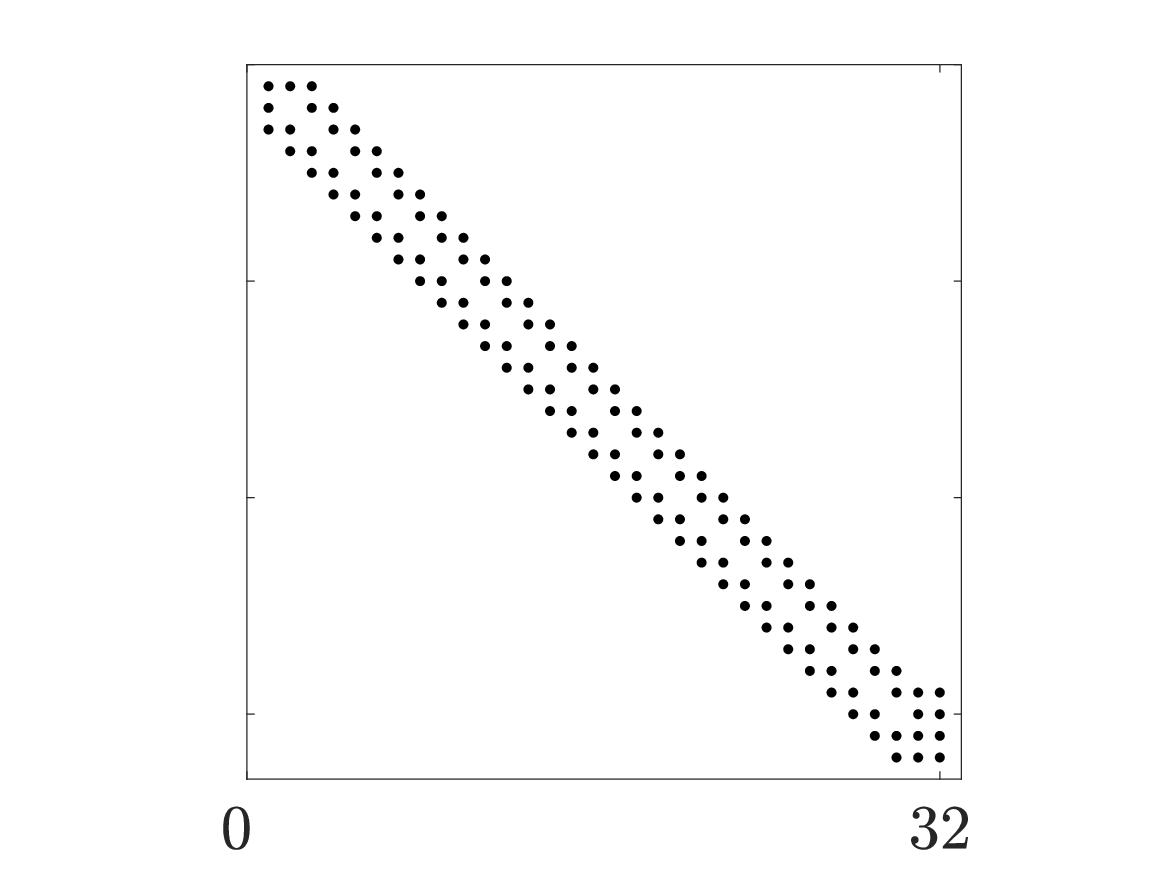}
  \caption{Sparsity pattern of the $B$ matrix.}
  \label{fig:bSparsity}
\end{subfigure}
\hfill
\begin{subfigure}[t]{50mm}
  \includegraphics[width=50mm]{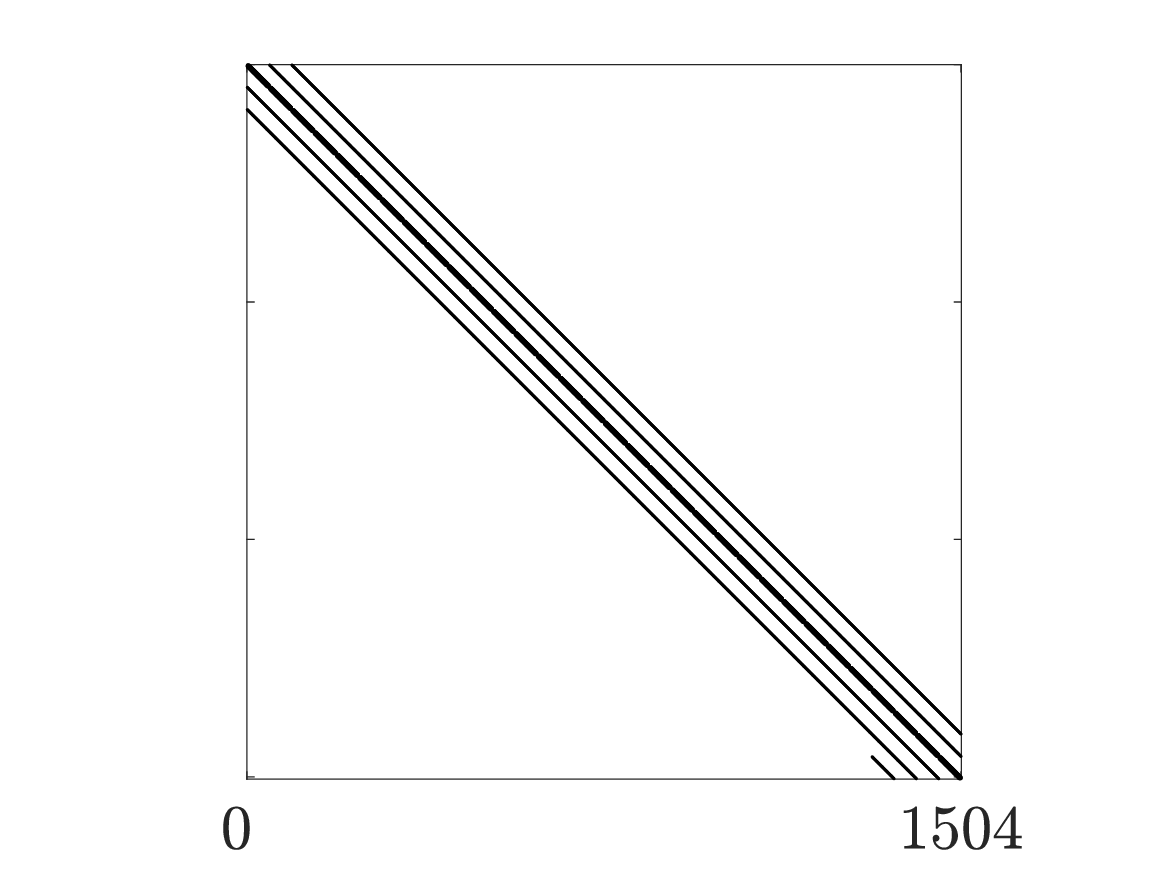}
  \caption{Sparsity pattern of the $K$ matrix.}
  \label{fig:kSparsity}
\end{subfigure}
\centering
\caption{Sparsity patterns of the $A$, $B$, (Eq. (\ref{eq:syl})) and $K$ (Eq. (\ref{eq:kron})) matrices obtained by discretizing the linear diffusion equation, Eq. \ref{eq:linDiff}.}
\label{fig:sparsity}
\end{figure}

To illustrate spacetime wavelet discretization, consider the 1D linear diffusion equation
\begin{align}
\begin{split}
    &\frac{\partial f(\vec{q})}{\partial t} - \nu \frac{\partial^2f(\vec{q}) }{\partial x^2} = 0,  \ \ \mathrm{in} \ \ \Omega, \\ 
    & f(\vec{q}) = f_I \ \ \mathrm{on} \ \ \Omega_x \times (t=0), \\
    &f(\vec{q}) = f_B \ \ \mathrm{on} \ \ \partial\Omega_x \times (0,T].
    \end{split}
    \label{eq:LD}
\end{align}
When discretized with dense wavelet derivative operators, the problem is expressed in matrix form as
\begin{align} \label{eq:linDiff}
    \boldsymbol{\mathcal{F}}\cdot {}^{(1,t)}\boldsymbol{\Gamma} - \nu {}^{(2,x)}\boldsymbol{\Gamma}\cdot \boldsymbol{\mathcal{F}} = \mathbf{0}.
\end{align}
Note that the ${}^{(2,x)}\boldsymbol{\Gamma}$ operator that approximates the spatial second derivative operates on the first index of $\boldsymbol{\mathcal{F}}$, corresponding to the spatial dimension, and ${}^{(1,t)}\boldsymbol{\Gamma}$ on the second, temporal index. The structure of Eq. (\ref{eq:linDiff}) is of the form 
\begin{align} 
    AX + XB = C,
    \label{eq:syl}
\end{align}
classified as a Sylvester equation \cite{bartels1972solution}, where  $A \in \mathbb{R}^{n\times n}$, $B \in \mathbb{R}^{s\times s}$, $C\in \mathbb{R}^{n\times s}$, $X \in \mathbb{R}^{n\times s}$. As mentioned previously, the Sylvester equation can be converted into a $K\vec{x}=\vec{r}$ problem using the Kronecker product of two matrices and matrix vectorization denoted by $(\bullet) \otimes (\bullet)$ and $\text{vec}(\bullet)$, respectively \cite{van2000ubiquitous}. This technique converts the problem to an $ns \times ns$ system
\begin{align}
    (I_n \otimes A + B^T \otimes I_s)\text{vec}(X) = \text{vec}(C),
    \label{eq:kron}
\end{align}
where $I_n$ is an $n\times n$ identity matrix. Because this equation is of the form $K\vec{x}=\vec{r}$, it can be solved with traditional linear solvers, given that $K$ is non-singular. Fig. \ref{fig:sparsity} shows the sparsity patterns for the $A$, $B$, and $K$ matrices resulting from the spacetime discretization of the linear diffusion equation, Eq. \ref{eq:linDiff}, in both the Sylvester form, Eq. (\ref{eq:syl}) and the Kronecker product form, Eq (\ref{eq:kron}), at resolution level $j=2$ with $p_x=6$, $p_t=4$. These wavelet parameters define the problem size: $n = 47$ and $s = 32$.
We see that all matrices are square and have a sparse banded structure. However the $K$ matrix shown in Fig. \ref{fig:kSparsity} is significantly larger than the $A$ and $B$ matrices (Figs. \ref{fig:aSparsity}, \ref{fig:bSparsity}) and contains more nonzero entries. The $A$ and $B$ matrices at $j=2$ contain 403 and 126 nonzeros, respectively,  while the $K$ matrix contains 18,677 nonzeros. This significant difference in number of nonzeros makes solving linear systems with the Kronecker product formulation computationally inefficient \cite{van2000ubiquitous}. Therefore, we proceed to solve the system in the original form of the Sylvester equation.

\subsection{Boundary Conditions} \label{sec:bcs}
The examples in this work are initial-boundary value problems with Dirichlet boundary conditions. Until the boundary and initial conditions are enforced, the problem is not well-posed and the $A$ and $B$ matrices (see Eq. (\ref{eq:syl})) are singular. We will enforce boundary condition constraints with a technique described by \cite{harker2015sylvester}, where semi-orthogonal permutation matrices separate the known and unknown degrees of freedom. We begin by defining the semi-orthogonal permutation matrices to operate on the spatial and temporal dimensions, $P_x$ and $P_t$, respectively, such that the generic solution $X$ can be expressed in terms of its known $X_D$ and unknown $\hat{X}$ components by
\begin{align}
    X = P_x \hat{X}P_t^T + X_D.
    \label{eq:bc1}
\end{align}
When substituted into the Sylvester form $AX+XB=C$, the equation becomes
\begin{align}
    AP_x\hat{X}P_t^T + P_x\hat{X}P_t^TB = C-AX_D - X_DB.
\end{align}
In order to leverage the semi-orthogonal properties of the permutation matrices to simplify the expression, we pre-multiply by $P_x^T$ and post-multiply by $P_t$ to obtain the equation
\begin{equation}
    \left(P_x^T  A  P_x\right)  \hat{X} + \hat{X} \left( P_t^T  B  P_t\right) = P_x^T(C-AX_D-X_DB)P_t.
\end{equation}
If we define 
\begin{align}
    \hat{A} &= P_x^TAP_x, & \hat{B} &= P_t^TBP_t, & \hat{C} &= P_x^T(C-AX_D-X_DB)P_t,
    \label{eq:bc2}
\end{align}
our equation is once again in the Sylvester form, $\hat{A}\hat{X}+\hat{X}\hat{B} = \hat{C}$, with square $\hat{A}$ and $\hat{B}$. The system is now well-posed with boundary and initial conditions enforced, and can be solved with compatible methods. An illustration of the $P_x$, $P_t$, $\hat{X}$ and $X_D$ matrices is presented in Appendix \ref{sec:bcMats}.

\subsection{Global-GMRES/Arnoldi} \label{sec:glGMRES}
The Gl-GMRES and Modified Global Arnoldi methods are used to solve the problems presented in this work and are described by Algorithms \ref{alg:glGMRES} and \ref{alg:mgArnoldi}, respectively. We denote the Frobenius inner product as $\langle \bullet , \bullet \rangle_F$ and the Frobenius norm as $\lVert \bullet \rVert_F$.
\begin{algorithm} [!htb]
\SetKwInOut{Input}{Inputs}
\SetKwInOut{Output}{Output}
\Input{$A$, $B$, $C$, $X_0$, $m$}
\Output{$X$}
Compute residual  $R_0 = C-AX_0-X_0B$. \\ 
Run Modified Global Arnoldi Algorithm  \Comment{Algorithm \ref{alg:mgArnoldi}} \\
Solve least squares problem for $\vec{y}$: \\
$\min\limits_{y\in \mathbb{R}^m} \| \|R_0\|_F \ \vec{e}_1 - H\vec{y} \|_2$, \ \ 
where $\vec{e}_1$ is the first unit vector, $e_1 = [1,0,0,...,0] \in \mathbb{R}^m$. \\ 
Update solution:  $X = X_0 + \mathcal{V}(\vec{y} \otimes I_s)$. \\ 
Set $X_0=X$.\\
Repeat until tolerance is met or maximum number of iterations is reached. \\
\caption{Global GMRES}
\label{alg:glGMRES}
\end{algorithm}
\begin{algorithm} [!htb]
\textbf{Inputs}: $A$, $B$, $C$, $X_0$, $R_0$, $m$, $\mathrm{tol}_H$ \\
\textbf{Outputs}: $\mathcal{V}$, $H$ \\
$V_1 = \frac{R_0}{\|R_0\|_F}$\\
\For{$z$ = 1\ :\ $m$}{
    $W = AV_z + V_z B$ \\ 
    \For{$i$ = 1 \ : \ $z$}{
        $H(i,z) = \langle V_i,W \rangle_F$ \\ 
        $W = W - H(i,z)V_i$ \\
    }
    $H(z+1,z) = \|W\|_F$ \\
    \If {$H(z+1,z) < \mathrm{tol}_H$}{
        break \\
    }
    $V_{z+1} = \frac{W}{H(z+1,z)}$ \\
}
$\mathcal{V} = [V_1, V_2, ... , V_m]$ \\ 
\caption{Modified Global Arnoldi Algorithm}
\label{alg:mgArnoldi}
\end{algorithm}
Algorithm \ref{alg:mgArnoldi} is referred to as ``modified" because it allows the user to specify the maximum number of orthogonal matrices computed, a parameter we denote as $m$. This is useful because, for many problems, the full span of the Krylov space is not required to obtain an accurate solution \cite{simoncini2000convergence, jaimoukha1997implicitly, beattie2005convergence}. Therefore, it is beneficial to keep the size of the Krylov space as small as possible while retaining enough orthogonal matrices to obtain an accurate solution, reducing memory requirements and computational cost. This feature makes the Gl-GMRES method a valuable tool when solving large Sylvester equations.

\subsection{Wavelet-Based Recursive Solution Technique}\label{sec:recursive}
Iterative solution techniques (\textit{e.g.}, Gl-GMRES) require an initial guess $X_0$ to begin the solution process. Previous research has shown that, for some iterative methods, a carefully-chosen initial guess can significantly speed up the computation time by reducing the number of iterations needed to reach a solution \cite{ye2020improving, elleithy2001iterative}. We propose a wavelet-based recursive algorithm to generate initial guesses and obtain these benefits. 

Our procedure begins by selecting a level $j<j_{\mathrm{max}}$ and solving the system using a generic initial guess (\textit{e.g.}, zeros). After obtaining the first solution, we use wavelet synthesis (see Eq. (\ref{eq:analSyn})) to construct the solution at the next level $j+1$. This solution serves as the initial guess for the next problem at level $j+1$. We repeat this process until we reach $j_{\mathrm{max}}$. This process results in fewer iterations required to reach a sufficiently accurate solution. The recursive wavelet algorithm is outlined by Algorithm \ref{alg:recAlg}.

\begin{algorithm}[H]
Read input\\
Define starting level $j$ and initial guess $X_0^j$\\
\While{$j \leq j_{\mathrm{max}}$}{
    Enforce boundary conditions \Comment{Eqs. (\ref{eq:bc1}-\ref{eq:bc2})} \\
    Solve Sylvester matrix equation using Gl-GMRES \Comment{Algorithms \ref{alg:glGMRES}, \ref{alg:mgArnoldi}} \\
    \If{$j < j_{\mathrm{max}}$}{Synthesize $X_0^{j+1}$ from solution \Comment{(Eq. \ref{eq:analSyn})}} 
    $j = j+1$
}
Compute error and output solution
\caption{Recursive Spacetime Wavelet Solver}
\label{alg:recAlg}
\end{algorithm} 

\section{Numerical Examples}
\label{sec:numerical}
In this section, we present the results of the numerical verification problems for the linear spacetime wavelet solver. The two problems we analyze are the dimensionless linear diffusion and convection-diffusion equations. For both problems, we use the Method of Manufactured Solutions (MMS) \cite{salari2000code}, allowing convergence analysis of the solution and derivative approximations. The reported errors for both examples are calculated by evaluating the highest wavelet coefficient on level $j_{\mathrm{max}} +1$, obtained with wavelet synthesis. This provides us with a reliable measure of the largest error on the grid at $j_{\mathrm{max}}$. We enforce boundary conditions and utilize the recursive wavelet algorithm to solve both problems using Gl-GMRES, as outlined in Sections \ref{sec:bcs} - \ref{sec:recursive}. We define the stopping criteria for the Gl-GMRES method as $\lVert R\rVert_F < 10^{-8}$ and the break criteria for the Modified Global Arnoldi algorithm as $\mathrm{tol}_H = 10^{-8}$. We use a restart value $m$ defined as a function of the resolution level $j$, \textit{i.e.}, $m = 30(j+1)$, to accommodate the increasing problem size. This value is chosen to achieve a balance of memory savings and solver performance. We denote the MMS and approximate wavelet solutions as $f_{ex}$ and $f_{\varepsilon}$, respectively. Both problems are solved on the domain $x \in[-1, 1]$, $t\in [0,1]$. All examples displayed in this work do not utilize preconditioning to present more equal comparisons.

\subsection{Linear Diffusion Equation} \label{sec:linearDiff}
The first example we present is the linear diffusion equation to test the capability of the spacetime wavelet solver to accurately simulate the diffusion of an arbitrary quantity (\textit{e.g.}, density), defined by Eq. (\ref{eq:LD}).
The boundary and initial conditions are defined by the analytical solution
\begin{align}
    f(x,t) = \sqrt{\frac{v^2\sigma^2}{2\nu t +\sigma^2}}\mathrm{exp}\left(\frac{-(x-x_0)^2}{2(2\nu t +\sigma^2)}\right),
\end{align}
 where $v = 1$, $\nu = 0.01$, $\sigma = 0.1$, and $x_0 = 0$. The parameter $v$ determines the amplitude of the solution, $\nu$ serves as a diffusion coefficient, $\sigma$ defines the width of the density peak and $x_0$ defines the center of the peak.

\begin{figure}[H]
\begin{subfigure}[t]{80mm}
  \includegraphics[width=80mm]{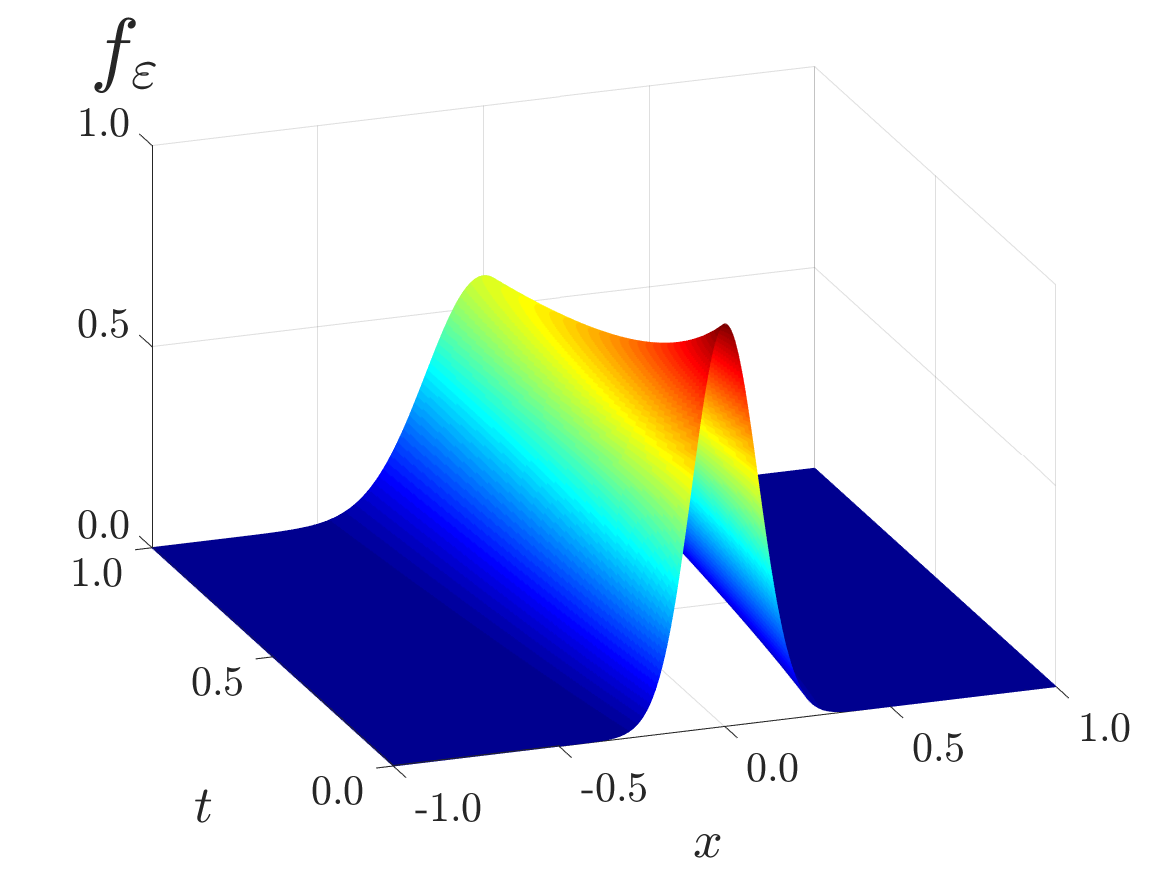}
  \caption{Spacetime solution.}
  \label{fig:ldST}
\end{subfigure}
\hfill
\begin{subfigure}[t]{80mm}
  \includegraphics[width=80mm]{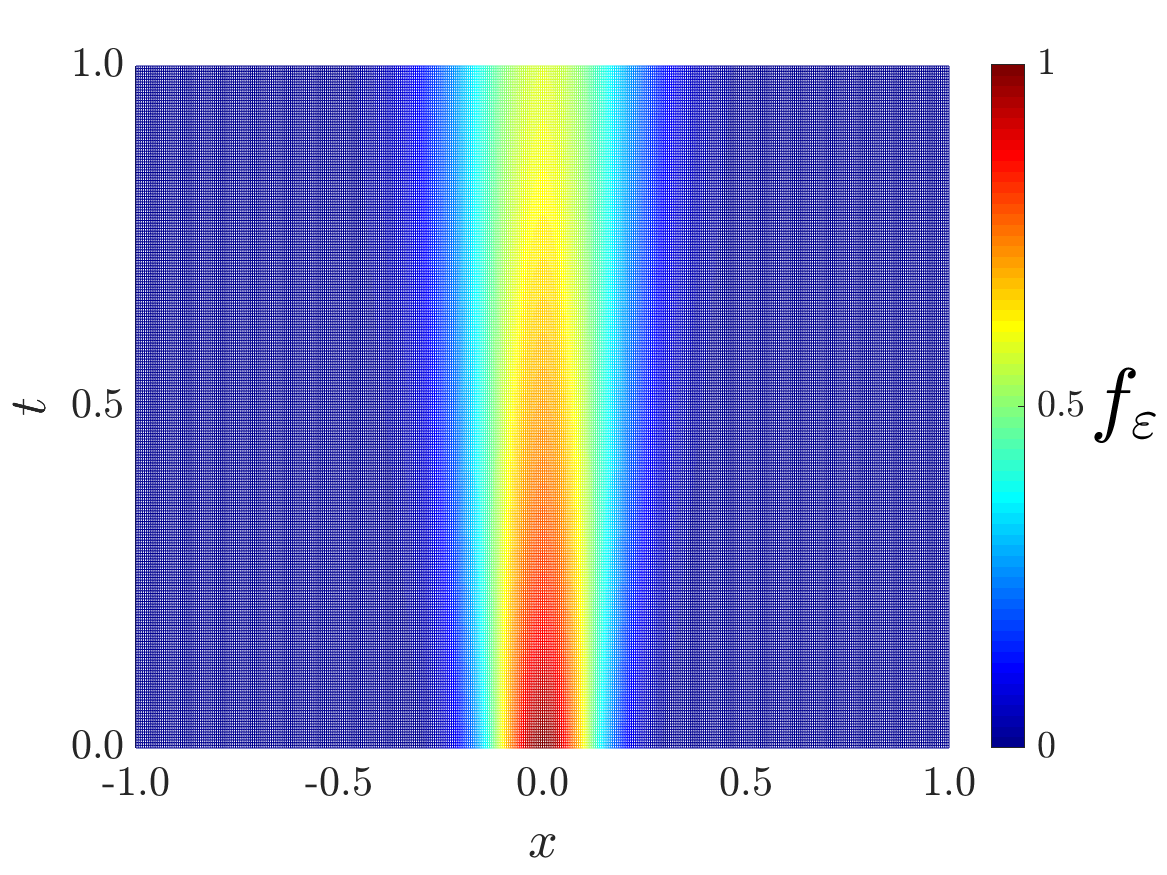}
  \caption{Top view.}
  \label{fig:ldTop}
\end{subfigure}
\centering
\caption{Spacetime solution for the linear diffusion problem at $j = 6$ with $p_x = 6$, $p_t = 4$.}
\label{fig:ld3d}
\end{figure}

 The spacetime solution at $j = 6$ is displayed in Fig. \ref{fig:ld3d}. We see that the initially steep peak present in the initial condition diffuses into a more shallow curve as the solution progresses through the time domain. This problem was solved with the basis order combination $p_x=6$, $p_t=4$. This combination is chosen to minimize the computational work required to build and solve the system while obtaining an accurate solution with high-order convergence.  An accuracy of $\lVert f_{ex}-f_{\varepsilon}\rVert_{\infty} = 6.1 \times 10^{-8}$ is obtained at $j=6$, solving for 392,703 degrees of freedom, denoted by $\mathcal{N}$. The $A$ and $B$ matrices contain 6,883 and 2,046 nonzero entries, respectively, with $1.05 \%$ combined sparsity. 

\begin{figure}[H]
\begin{subfigure}[t]{80mm}
  \includegraphics[width=80mm]{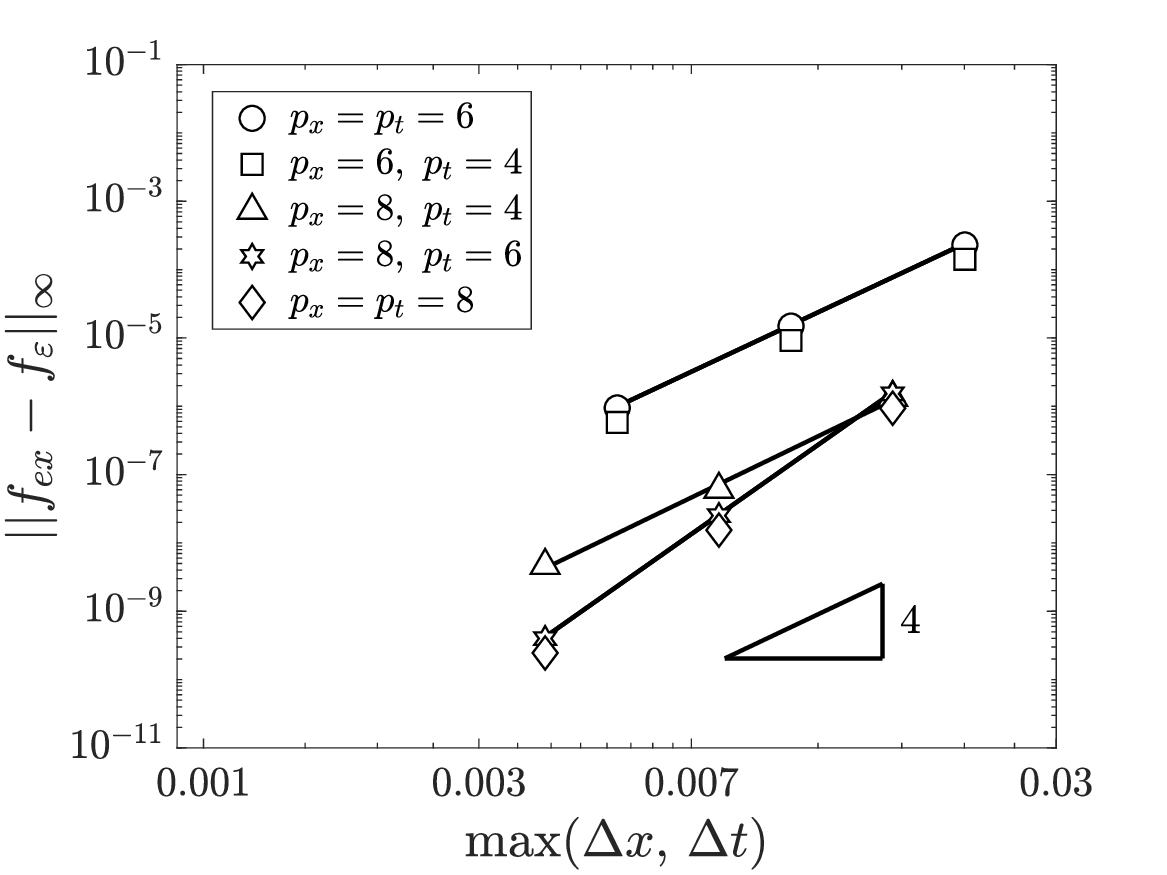}
  \caption{Solution convergence with varying $p_x$, $p_t$ combinations.}
  \label{fig:ldSolConv}
\end{subfigure}
\hfill
\begin{subfigure}[t]{80mm}
  \includegraphics[width=80mm]{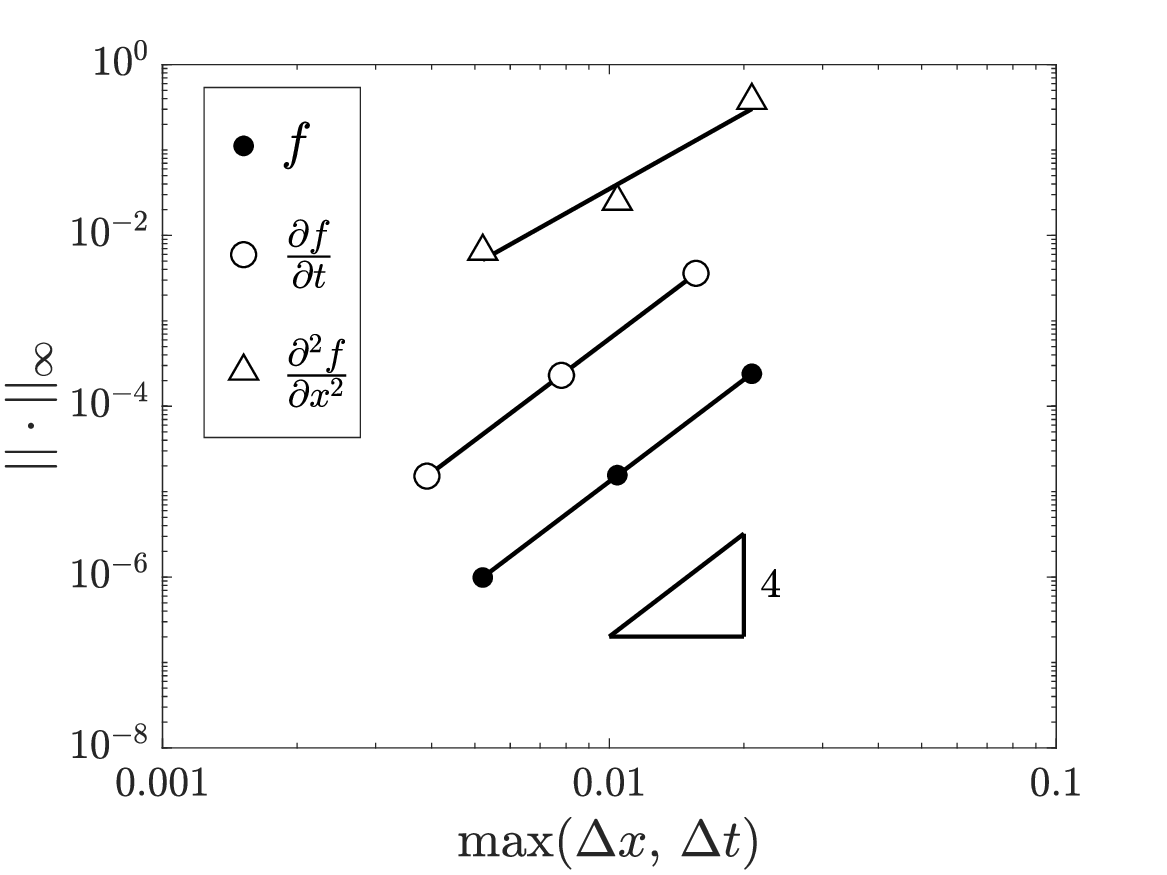}
  \caption{Solution and derivative convergence with $p_x=6$, $p_t=4$.}
  \label{fig:ldDconv}
\end{subfigure}
\centering
\caption{Solution and derivative convergence for the linear diffusion problem for $j = 3,4,5$.}
\label{fig:ldConv}
\end{figure}

Fig. \ref{fig:ldSolConv} displays the solution convergence rates for various combinations of spatial and temporal basis functions and Fig. \ref{fig:ldDconv} shows the convergence rates of the solution and derivative approximations at $p_x=6$, $p_t=4$. As one can see, when we increase the order of both the temporal and spatial basis functions, we obtain more accurate solutions that obey the predicted convergence rate of the function \cite{harnish2018adaptive, harnish2023adaptive}. Superconvergent behavior of the derivatives is illustrated in \ref{fig:ldDconv} with the derivative approximations converging at a rate higher than the expected order of 3, Eq. (\ref{eq:derErr}). We display the convergence rates for the linear diffusion solutions in Fig. \ref{fig:ldSolConv} in Table \ref{tab:ldConv}. 
\begin{table}[H]
\centering
\begin{tabular}{ |c|c|c| }
    \hline
    $p_x$ & $p_t$ & Convergence rate \\ 
    \hline
    6 & 6 & 3.96\\
    \hline
    6 & 4 & 3.96\\ 
    \hline
    8 & 4 & 4.08\\ 
    \hline
    8 & 6 & 5.95\\ 
    \hline
    8 & 8 & 5.95\\
    \hline
\end{tabular}
\captionsetup{justification=centering}
\caption{Convergence rates for the linear diffusion solution with $j = 3, 4, 5$ and varying $p_x$, $p_t$ (Fig. \ref{fig:ldSolConv}).}
\label{tab:ldConv}
\end{table}

\begin{figure}[H]
\begin{subfigure}[t]{50mm}
  \includegraphics[width=50mm]{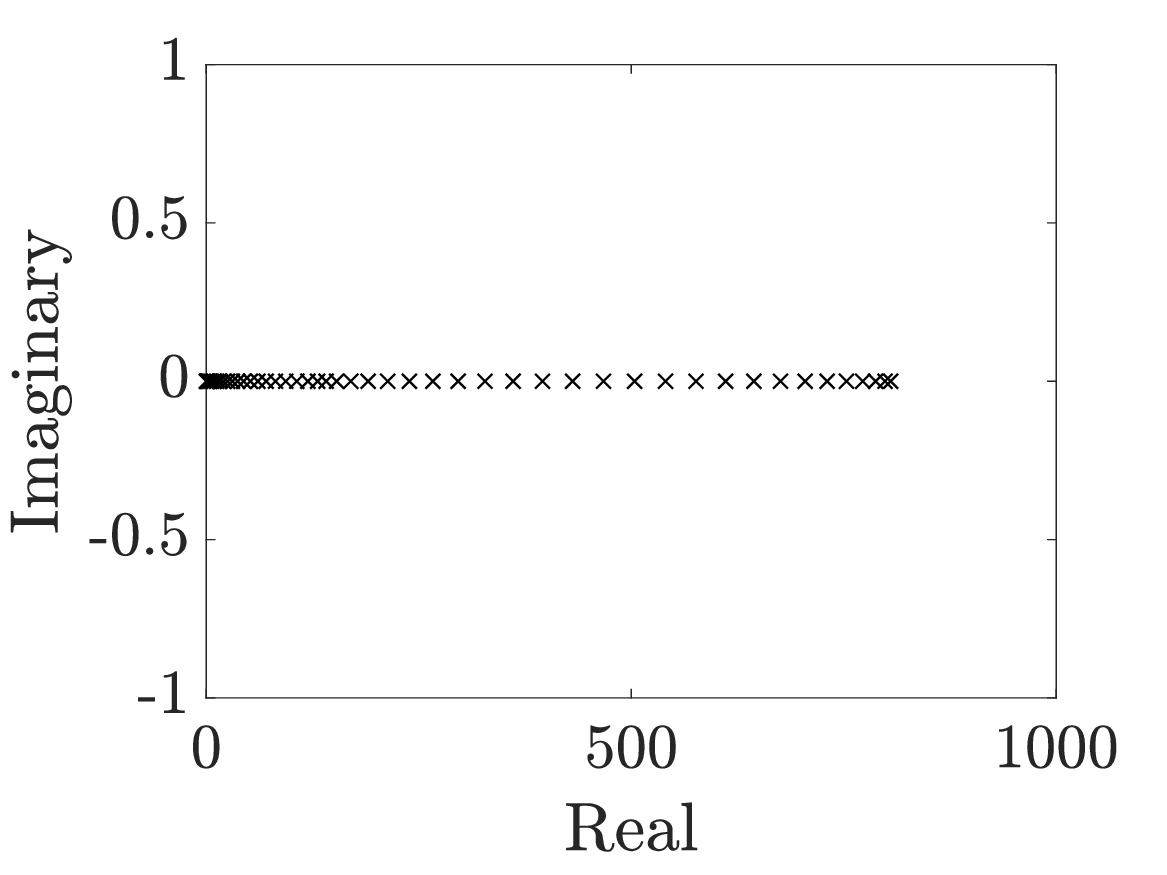}
  \caption{Eigenvalue spectrum of the $A$ matrix, Eq. (\ref{eq:syl}).}
  \label{fig:aSpec_ld}
\end{subfigure}
\hfill
\begin{subfigure}[t]{50mm}
  \includegraphics[width=50mm]{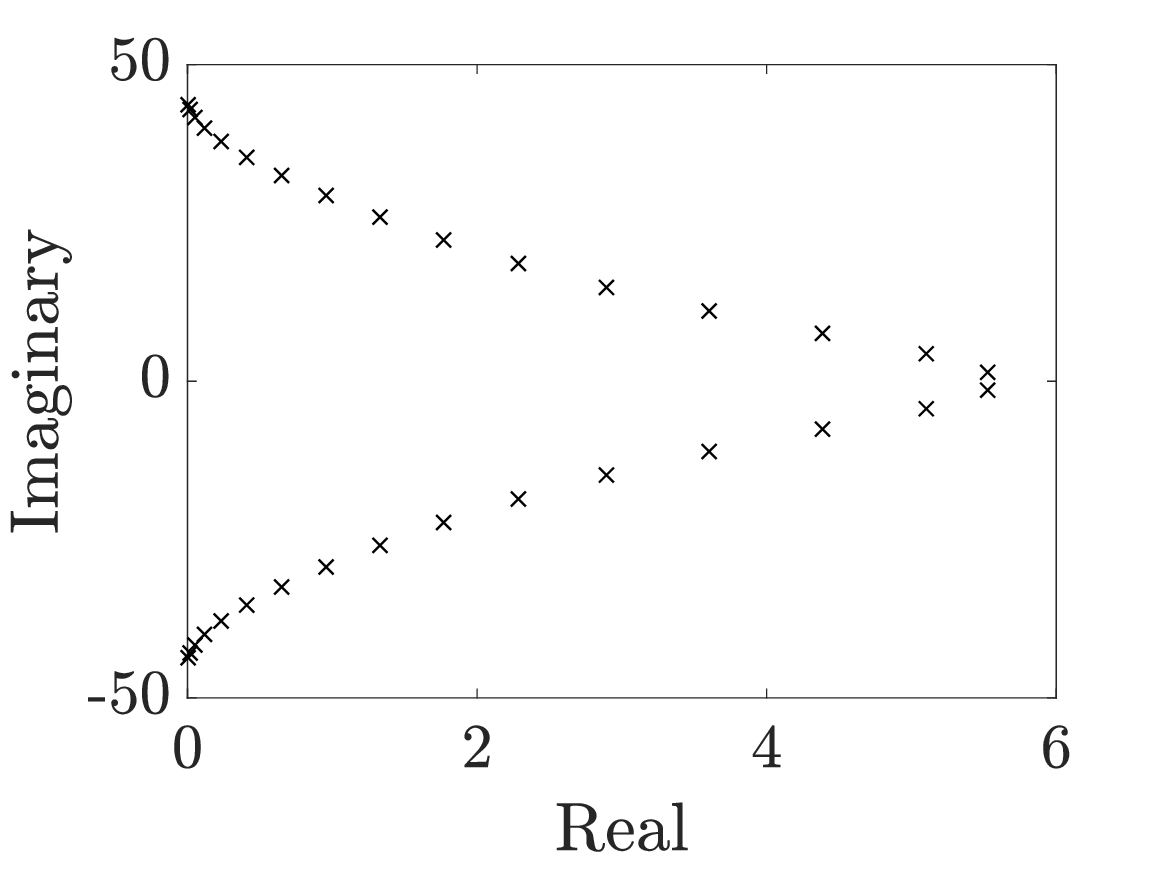}
  \caption{Eigenvalue spectrum of the $B$ matrix, Eq. (\ref{eq:syl}).}
  \label{fig:bSpec_ld}
\end{subfigure}
\hfill
\begin{subfigure}[t]{50mm}
  \includegraphics[width=50mm]{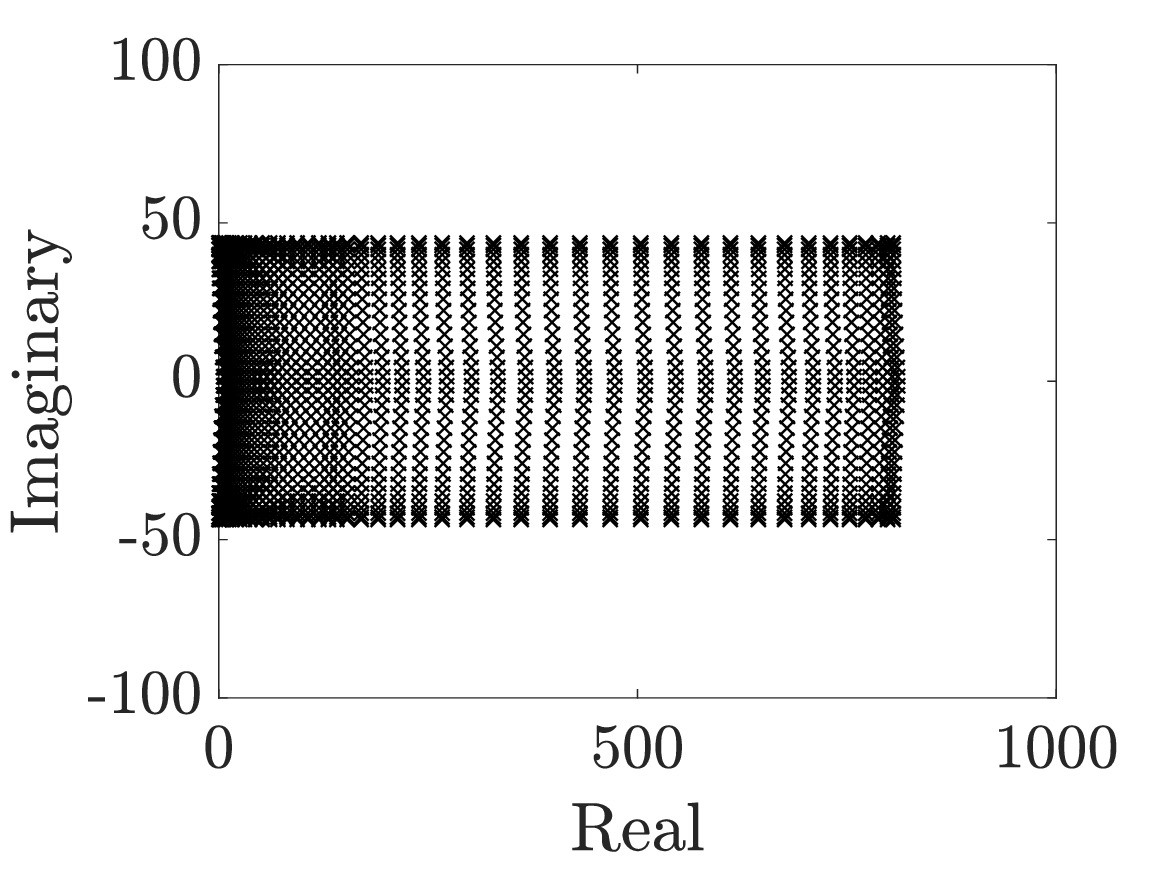}
  \caption{Eigenvalue spectrum of the $K$ matrix, Eq. (\ref{eq:kron}).}
  \label{fig:kSpec_ld}
\end{subfigure}
\centering
\caption{Eigenvalue spectra of the $A$, $B$, and $K$ matrices obtained by discretizing the linear diffusion equation, Eq. (\ref{eq:linDiff}).}
\label{fig:ldSpec}
\end{figure}

Fig. \ref{fig:ldSpec} shows the eigenvalue spectra of the $A$, $B$, and $K$ matrices.
For this example, the matrix $A$ contains only real eigenvalues. We also note that all real components of the eigenvalues for this problem are positive, opposite to what is typically observed in stability analysis of semi-discretized PDEs \cite{zauderer2011partial}. This occurs due to the use of the residual form of the governing equations, Eq. (\ref{eq:linDiff}), obtained from the spacetime discretization. We see that both formulations result in matrices with complex eigenvalues, illustrating why it is crucial to carefully select a compatible solver for both formulations \cite{saad2003iterative}. Detailed discussion about the eigenvalues and their dependency on the basis orders can be found in \cite{cochran2025spacetime}.

\subsection{Linear Convection-Diffusion Equation}
\begin{figure}[H]
\begin{subfigure}[t]{80mm}
  \includegraphics[width=80mm]{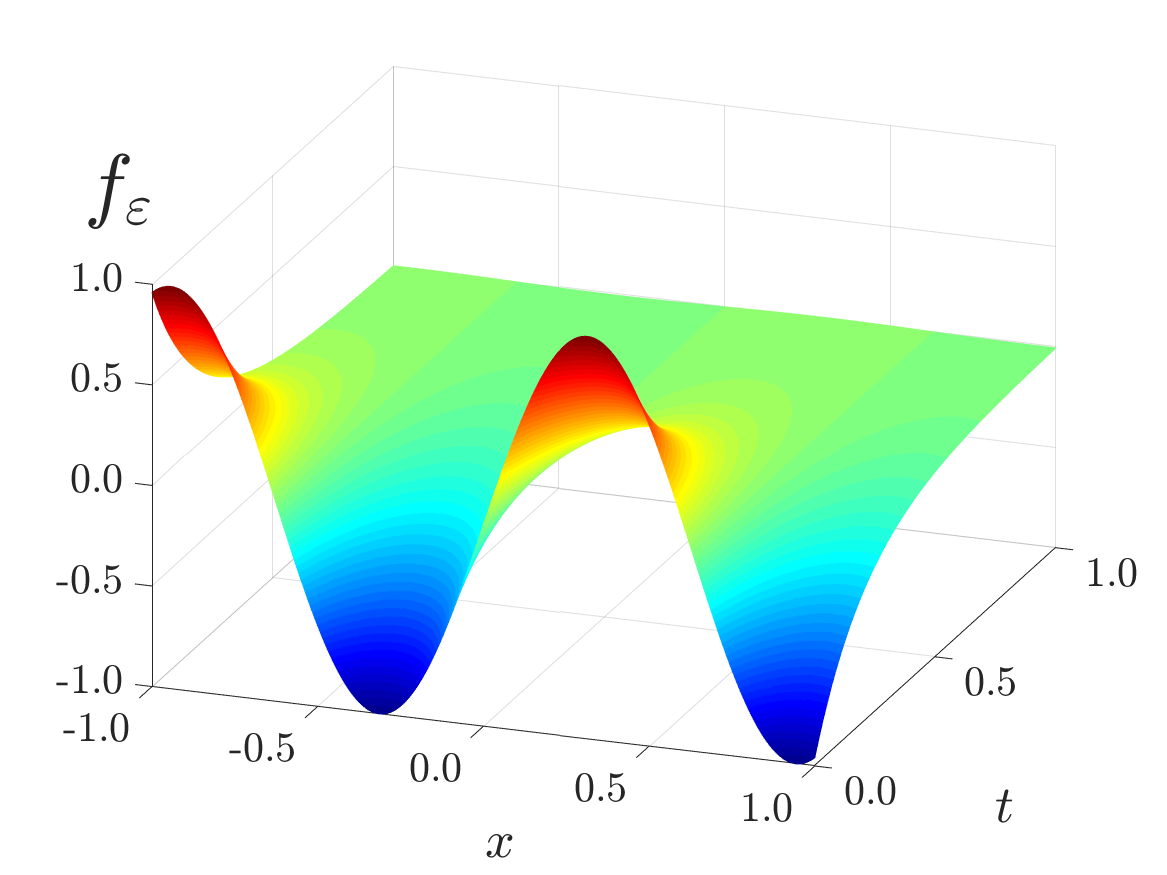}
  \caption{Spacetime solution.}
  \label{fig:cdST}
\end{subfigure}
\hfill
\begin{subfigure}[t]{80mm}
  \includegraphics[width=80mm]{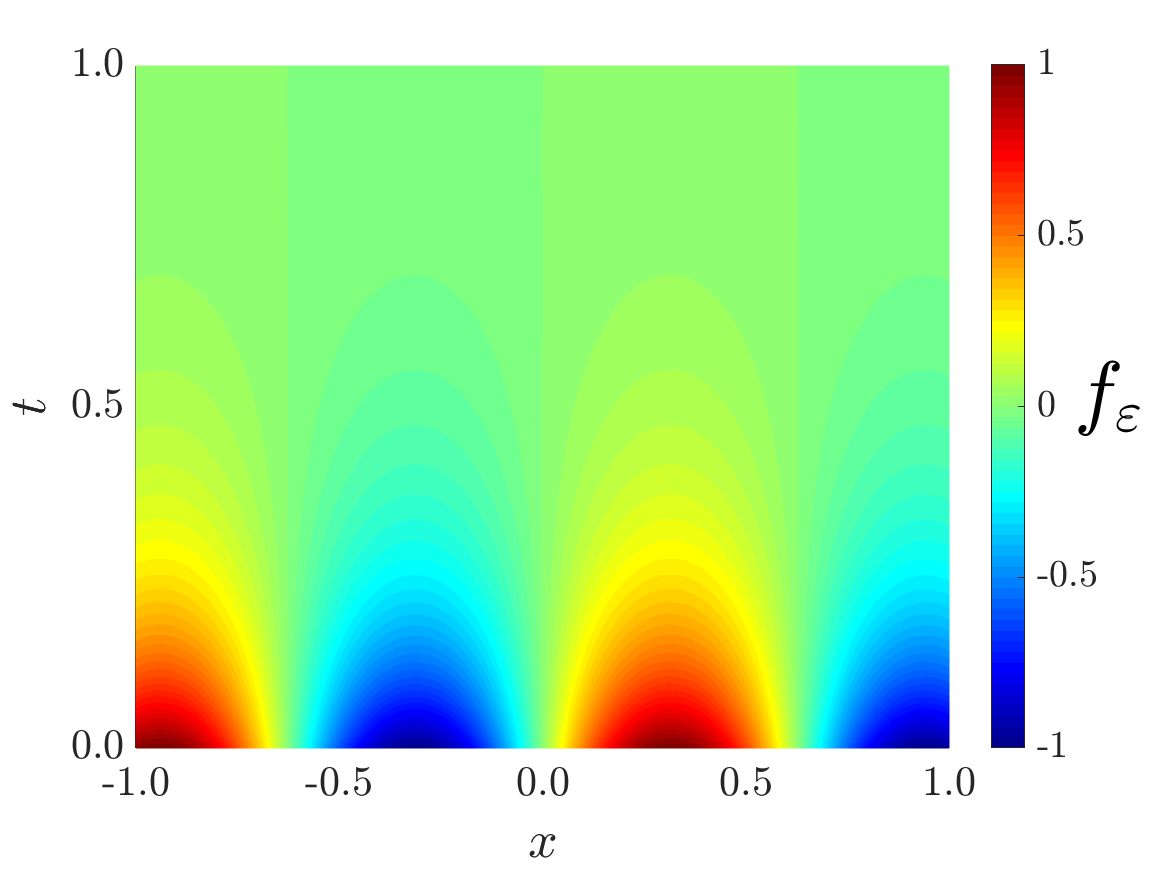}
  \caption{Top view.}
  \label{fig:cdTop}
\end{subfigure}
\centering
\caption{Spacetime solution for the convection-diffusion problem at $j = 6$ with $p_x = p_t=8$.}
\label{fig:cd3d}
\end{figure}

To demonstrate the capability of the spacetime wavelet solver to handle the addition of spatial convective terms, we discretize and solve the convection-diffusion equation. The convection-diffusion equation is governed by
\begin{align}
\begin{split}
    &\frac{\partial f(\vec{q})}{\partial t} + c\frac{\partial f(\vec{q})}{\partial x} - \nu \frac{\partial^2f(\vec{q}) }{\partial x^2} = 0,  \ \ \mathrm{in} \ \ \Omega, \\ 
    & f(\vec{q}) = f_I \ \ \mathrm{on} \ \ \Omega_x \times (t=0), \\
    &f(\vec{q}) = f_B \ \ \mathrm{on} \ \ \partial\Omega_x \times (0,T],
    \end{split}
\end{align}
with the boundary and initial conditions defined by the manufactured solution
\begin{align}
    f(x,t) = v \ \text{sin}\left(ax\right)\mathrm{e}^{-bt}.
    \label{eq:cdExSol}
\end{align}
\noindent For this problem, we define the problem parameters as $v = 3$, $a=5$, $b =5$, $c =  1$ , $\nu = 0.01$. From our selected analytical solution, Eq. (\ref{eq:cdExSol}), MMS yields the forcing term
\begin{align}
    g(x,t) = v\left[ac \ \mathrm{cos}(ax) + (a^2\nu - b) \ \mathrm{sin}(ax)\right]\mathrm{e}^{-bt}.
\end{align}
The discretized problem is expressed in matrix form as
\begin{align} \label{eq:conDiff}
        \boldsymbol{\mathcal{F}}\cdot{}^{(1,t)}\boldsymbol{\Gamma} + \left(c{}^{(1,x)}\boldsymbol{\Gamma} - \nu {}^{(2,x)}\boldsymbol{\Gamma}\right)\cdot \boldsymbol{\mathcal{F}} = \mathbf{G}.
\end{align}
For this example, we use the basis orders $p_x =p_t=8$ to show how increased basis orders affect accuracy and convergence. Fig. \ref{fig:cd3d} shows the spacetime solution of the convection-diffusion problem, which achieves an accuracy of $\lVert f_{ex} - f_{\varepsilon}\rVert_{\infty} = 1.96\times 10^{-9}$ at $j=6$ ($\mathcal{N} =$ 1,047,552). 

\begin{figure}[H]
  \begin{subfigure}[t]{75mm}
        \includegraphics[width=75mm]{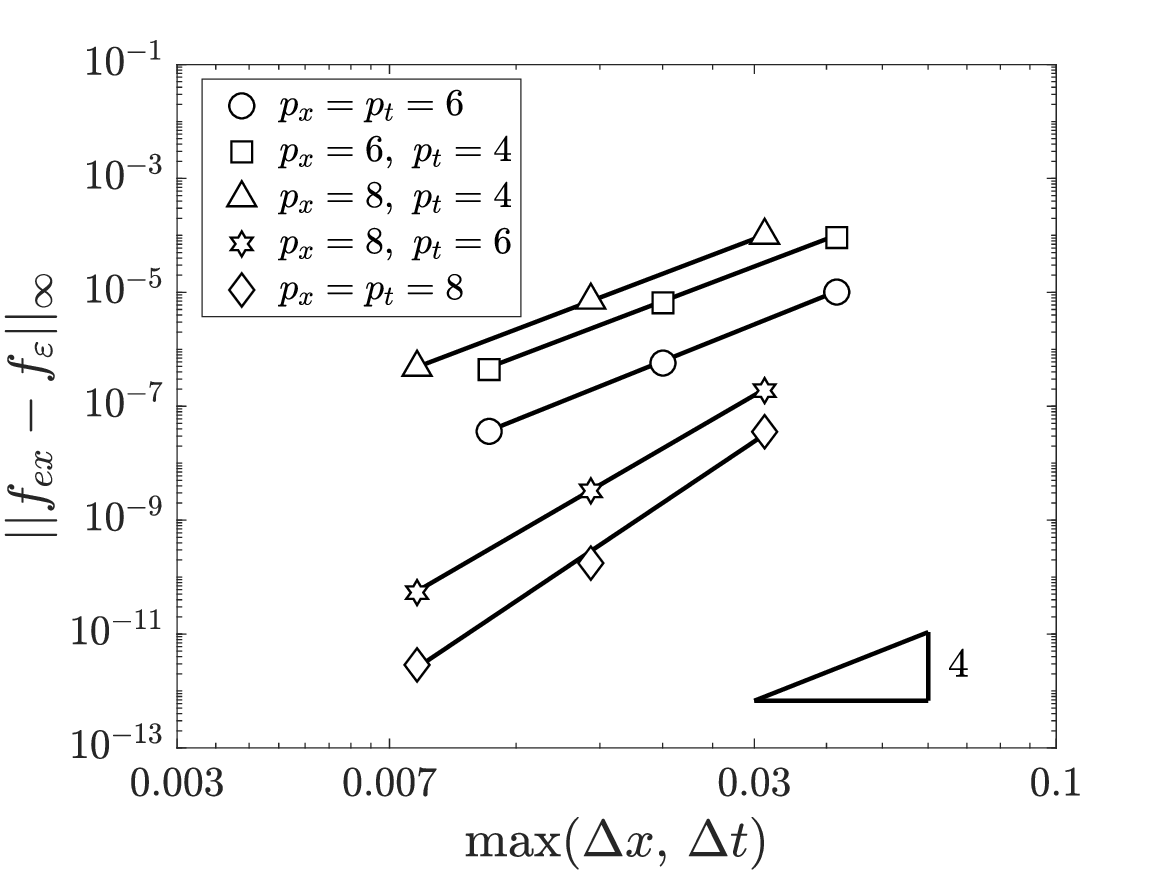}
        \caption{Solution convergence with varying $p_x$, $p_t$ combinations.}
        \label{fig:cdSolConv}
  \end{subfigure}
  \hfill
  \begin{subfigure}[t]{75mm}
        \includegraphics[width=75mm]{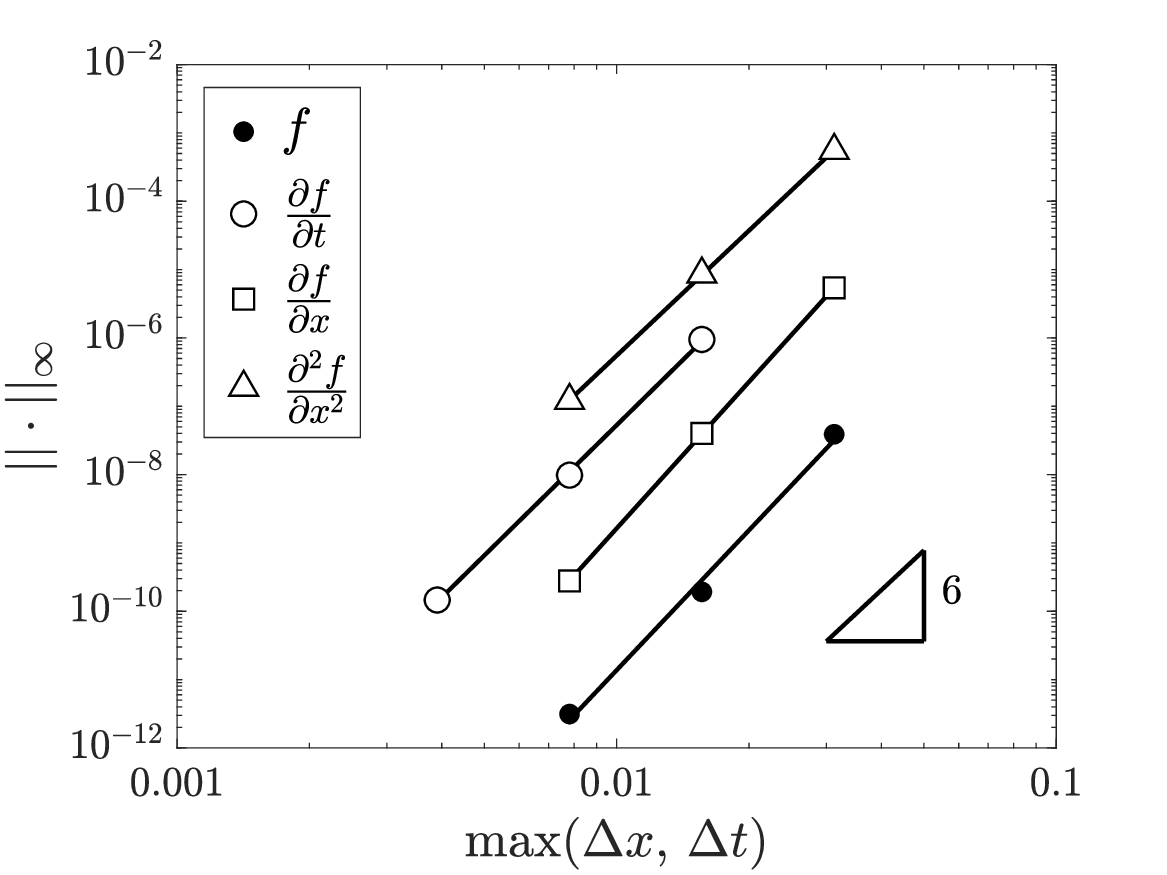}
        \caption{Solution and derivative convergence with $p_x=p_t=8$.}
        \label{fig:cdDconv}
  \end{subfigure}
  \centering
  \caption{Solution and derivative convergence for the convection-diffusion problem for $j = 2,3,4$.}
  \label{fig:cdConv}
\end{figure}

\begin{figure}[H]
\begin{subfigure}[t]{50mm}
  \includegraphics[width=50mm]{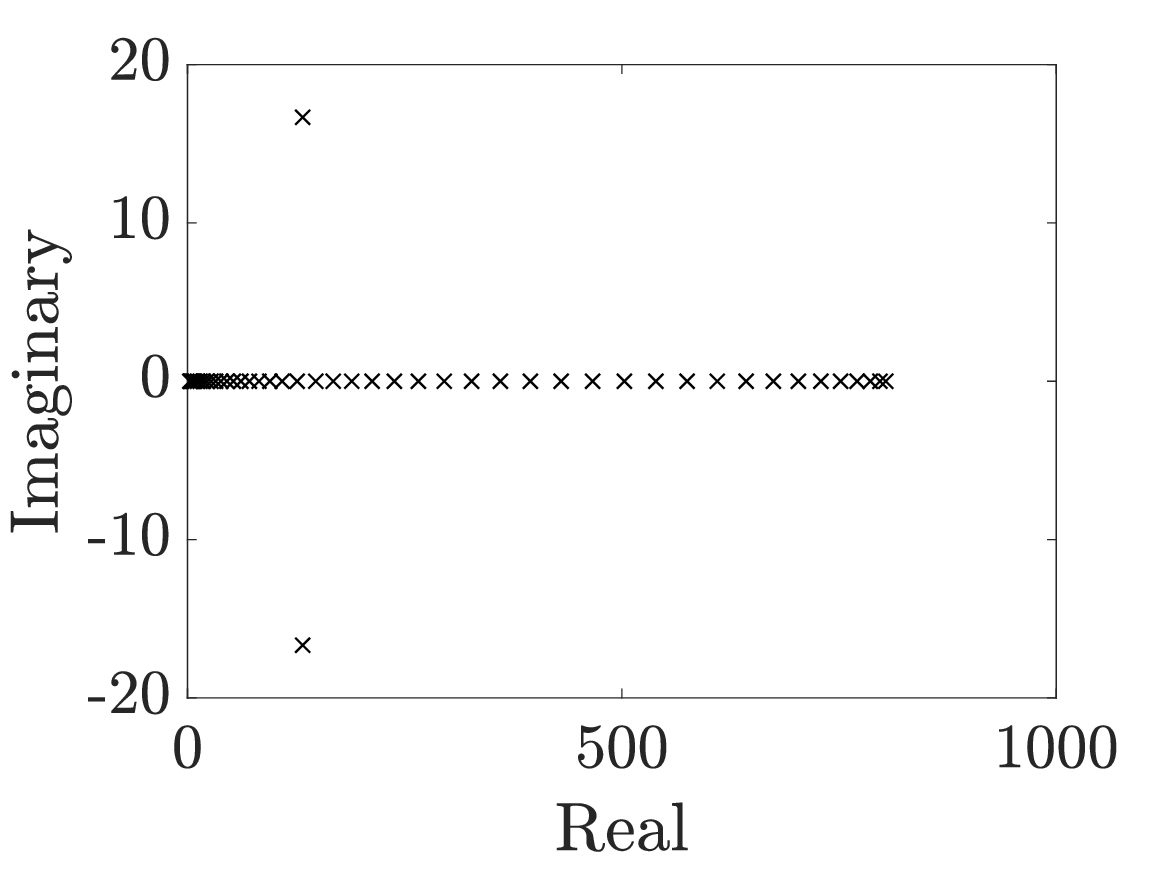}
  \caption{Eigenvalue spectrum of the $A$ matrix, Eq. (\ref{eq:syl}).}
  \label{fig:aSpec_cd}
\end{subfigure}
\hfill
\begin{subfigure}[t]{50mm}
  \includegraphics[width=50mm]{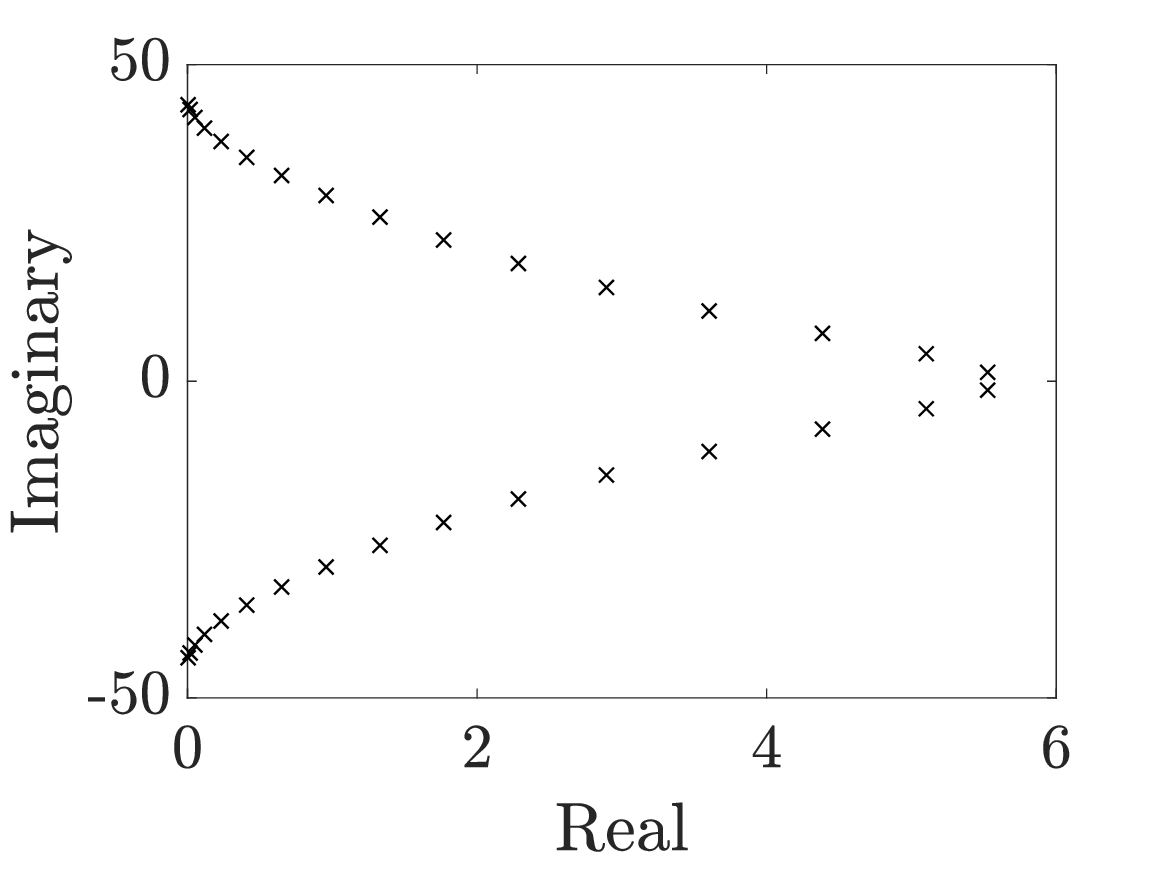}
  \caption{Eigenvalue spectrum of the $B$ matrix, Eq. (\ref{eq:syl}).}
  \label{fig:bSpec_cd}
\end{subfigure}
\hfill
\begin{subfigure}[t]{50mm}
  \includegraphics[width=50mm]{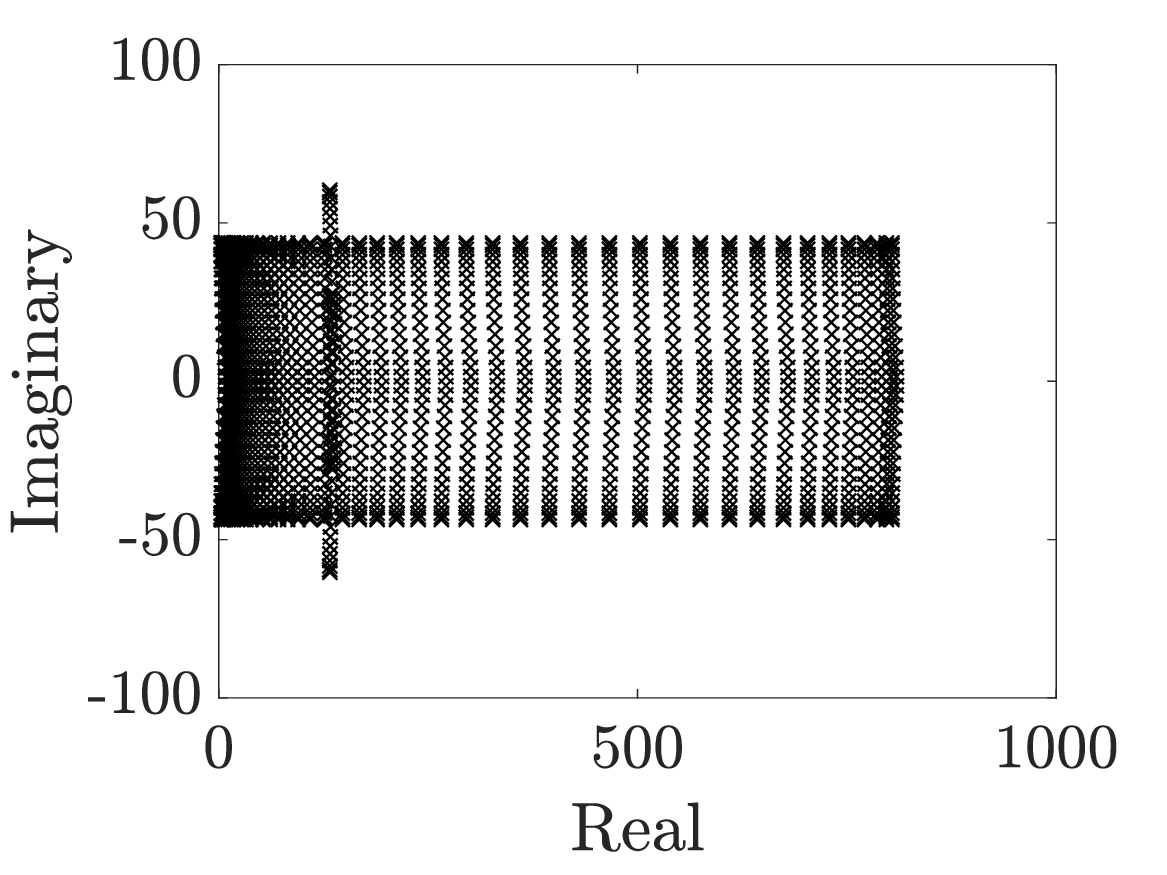}
  \caption{Eigenvalue spectrum of the $K$ matrix, Eq. (\ref{eq:kron}).}
  \label{fig:kSpec_cd}
\end{subfigure}
\centering
\caption{Eigenvalue spectra of the $A$, $B$, and $K$ matrices obtained by discretizing the convection-diffusion equation (Eq. \ref{eq:conDiff}) at $j=2$ with $p_x=6$, $p_t=4$.}
\label{fig:cdSpec}
\end{figure}

As with the linear diffusion example, we display solution and convergence rates in Fig. \ref{fig:cdConv}. As predicted by the theory, Eq. (\ref{eq:derErr}), we achieve $6$th-order convergence or better for both the solution and derivative approximations. Note that for Fig. \ref{fig:cdConv}, we show the results for $j=2,3,4$.
We display the convergence rates for the convection-diffusion simulations in Table \ref{tab:cdConv}.
\begin{table}[H]
\centering
\begin{tabular}{ |c|c|c| }
    \hline
    $p_x$ & $p_t$ & Convergence rate \\ 
    \hline
    6 & 6 & 4.06\\
    \hline
    6 & 4 & 3.85\\ 
    \hline
    8 & 4 & 3.85\\ 
    \hline
    8 & 6 & 5.89\\ 
    \hline
    8 & 8 & 6.80\\
    \hline
\end{tabular}
\captionsetup{justification=centering}
\caption{Convergence rates for the convection-diffusion solution with $j = 3, 4, 5$ and varying $p_x$, $p_t$ (Fig. \ref{fig:cdSolConv}).}
\label{tab:cdConv}
\end{table}
We see similar trends with this example as with the linear diffusion example, where increasing basis orders increases convergence rates. Superconvergence is again obtained for the derivatives. Fig. \ref{fig:cdSpec} shows the eigenvalue spectra of the $A$, $B$, and $K$ matrices for the convection-diffusion problem at $j=2$ with $p_x=p_t=8$.
We see that the addition of the convection term results in complex eigenvalues in the $A$ matrix (Fig. \ref{fig:aSpec_cd}) and the appearance of a vertical band seen in the $K$ matrix (Fig. \ref{fig:kSpec_cd}), not seen with the linear diffusion example (Fig. \ref{fig:ldSpec}). These differences in the spectra are present for all combinations of $p_x$ and $p_t$.

This section has demonstrated the capability of the spacetime wavelet solver to effectively and accurately discretize and solve linear PDEs in both the spatial and temporal dimensions simultaneously. The high-order convergence predicted by the wavelet theory, see Eq. (\ref{eq:derErr}), is achieved for both the solution and derivative approximations.

\subsection{Comparison of Kronecker product formulation and the Sylvester Equation} \label{sec:axbVaxxbc}
\begin{figure}[H]
\begin{subfigure}[t]{75mm}
  \includegraphics[width=75mm]{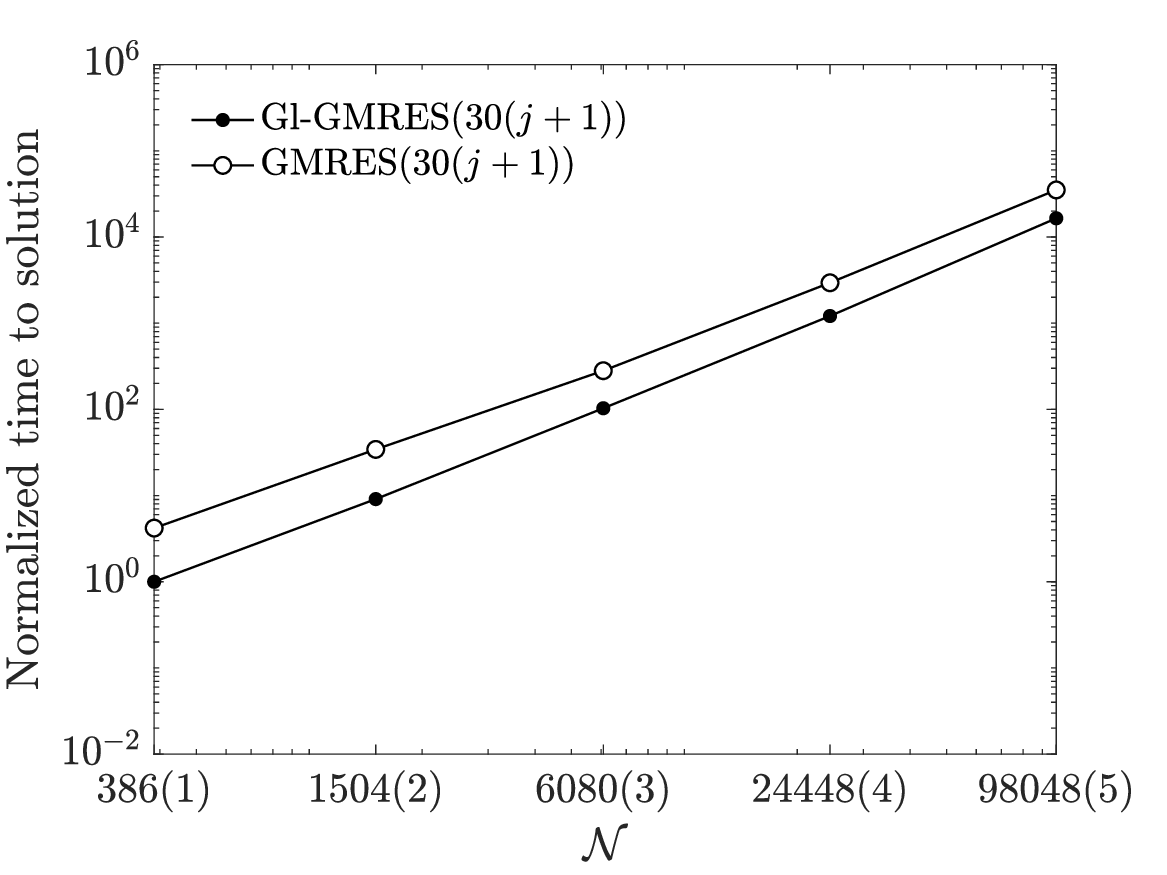}
  \caption{Time to solution normalized with the $j=1$ Gl-GMRES time.}
  \label{fig:ldaxb30}
\end{subfigure}
\hfill
\begin{subfigure}[t]{75mm}
  \includegraphics[width=75mm]{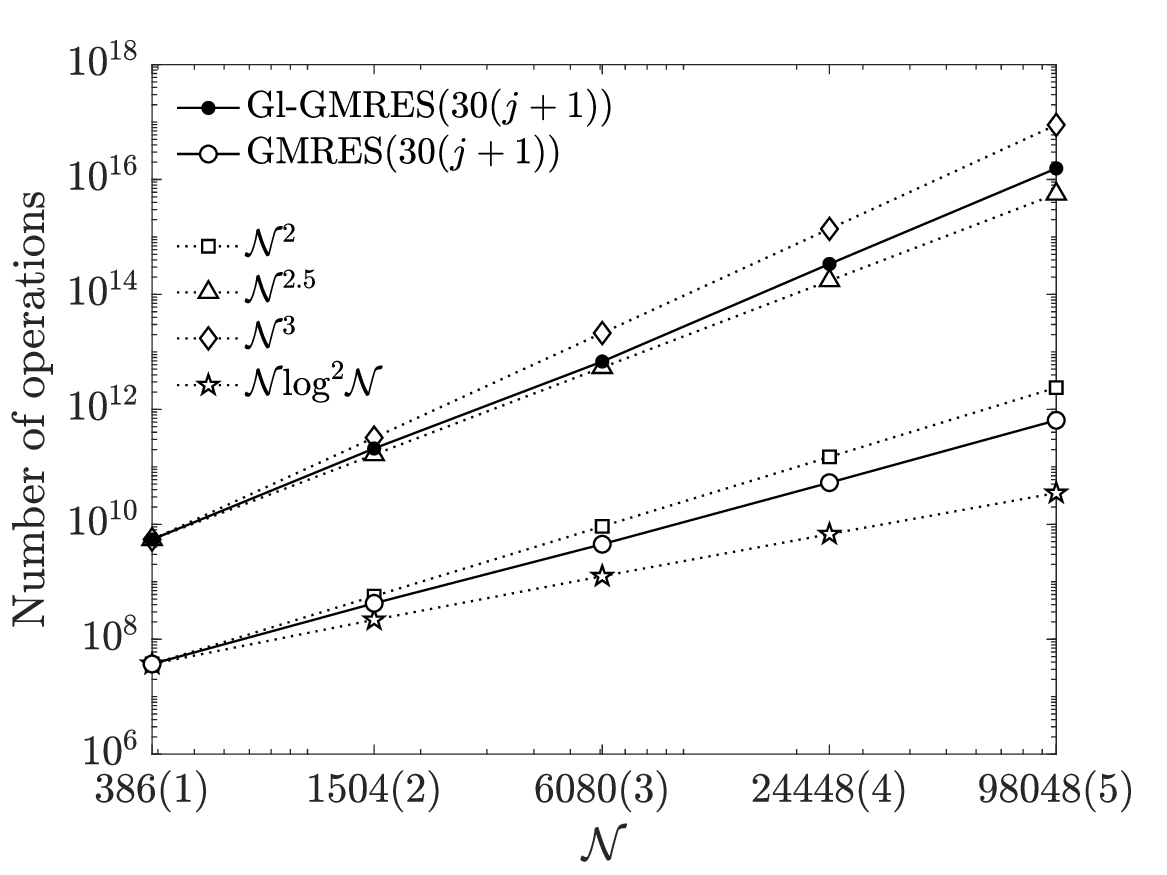}
  \caption{Number of operations as a function of $\mathcal{N}$.}
  \label{fig:nops}
\end{subfigure}
\centering
\caption{Performance comparisons of the Gl-GMRES and standard restarted GMRES methods using restart parameter $m = 30(j+1)$, $p_x=6$, $p_t=4$ at levels $j=1-5$. Level $j$ is displayed in parentheses.}
\label{fig:axb}
\end{figure}

In this section, we provide a comparison between two techniques used to solve the Sylvester matrix equation resulting from our linear diffusion problem in Section \ref{sec:linearDiff}. We compare the time required to assemble and solve the system using the Kronecker product formulation, Eq. (\ref{eq:kron}), with a standard restarted GMRES method \cite{simoncini2000convergence} to the Sylvester equation solved with the Gl-GMRES method \cite{jbilou1999global}. For a fair comparison, we use an initial guess of zeros for both methods and set the maximum size of the Krylov space to the same value $m=30(j+1)$. We solve with the advection parameter $c = 1$ and relaxed viscosity parameter $\nu=0.1$ to obtain accurate solutions at smaller $j$. Fig. \ref{fig:ldaxb30} shows the normalized time to solution for the linear diffusion example with times presented as an average over multiple runs from $j=1$ ($\mathcal{N}=$ 368) to $j=5$ ($\mathcal{N}=$ 98,048), using $p_x=6$, $p_t=4$. All times are normalized with the Gl-GMRES $j=1$ time. As one can see, Gl-GMRES solves the Sylvester equation consistently faster than GMRES solves $K\vec{x} = \vec{r}$. Fig. \ref{fig:nops} shows the number of operations required for both techniques. We see that the Gl-GMRES is close to $\mathcal{O}(\mathcal{N}^{2.5})$ and the restarted GMRES is slightly below $\mathcal{O}(\mathcal{N}^{2})$.  The restarted nature of GMRES variants significantly reduces memory requirements. However, they do not have a fixed order of convergence \cite{lin2005implicitly, jbilou1999global}, and therefore do not have a simply-determined algorithmic complexity. The number of iterations and total number of operations required for a convergent solution is highly dependent on the restart parameter $m$, characteristics of the matrices, initial guess, and residual tolerance \cite{baker2005technique, embree2003tortoise}. To put our computational complexity into context with other well-known numerical methods, the band Cholesky scheme is $\mathcal{O}(\mathcal{N}^{2})$, the CG method is $\mathcal{O}(\mathcal{N}^{1.5})$, and the preconditioned CG method is $\mathcal{O}(\mathcal{N}^{1.25})$ \cite{heath2018scientific}. Note that neither of the methods used in this work are preconditioned, and we expect that reduced computational cost would be achieved by implementing carefully-chosen preconditioning \cite{bouhamidi2008note, bouhamidi2013preconditioned, chen2005matrix}.

For all levels of the linear diffusion problem, building and solving the Sylvester equation with Gl-GMRES is faster than assembling the matrix, vectorizing, and solving the $K\vec{x}=\vec{r}$ problem with a restarted GMRES technique. At level $j=5$, the Sylvester equation $A$ and $B$ matrices require the storage of 3,427 and 1,022 nonzeros, respectively. The $K$ matrix from the Kronecker product formulation contains 1,267,589 nonzeros, three orders of magnitude more than the analogous Sylvester equation. This massive difference in required memory alone illustrates an advantage of solving large systems in the Sylvester form. We obtain similar results to those shown in Fig. \ref{fig:axb} when using different values of $m$.

\subsection{Recursive Solution Technique} \label{sec:recSol}
\begin{figure}[H]
\begin{subfigure}[t]{75mm}
  \includegraphics[width=75mm]{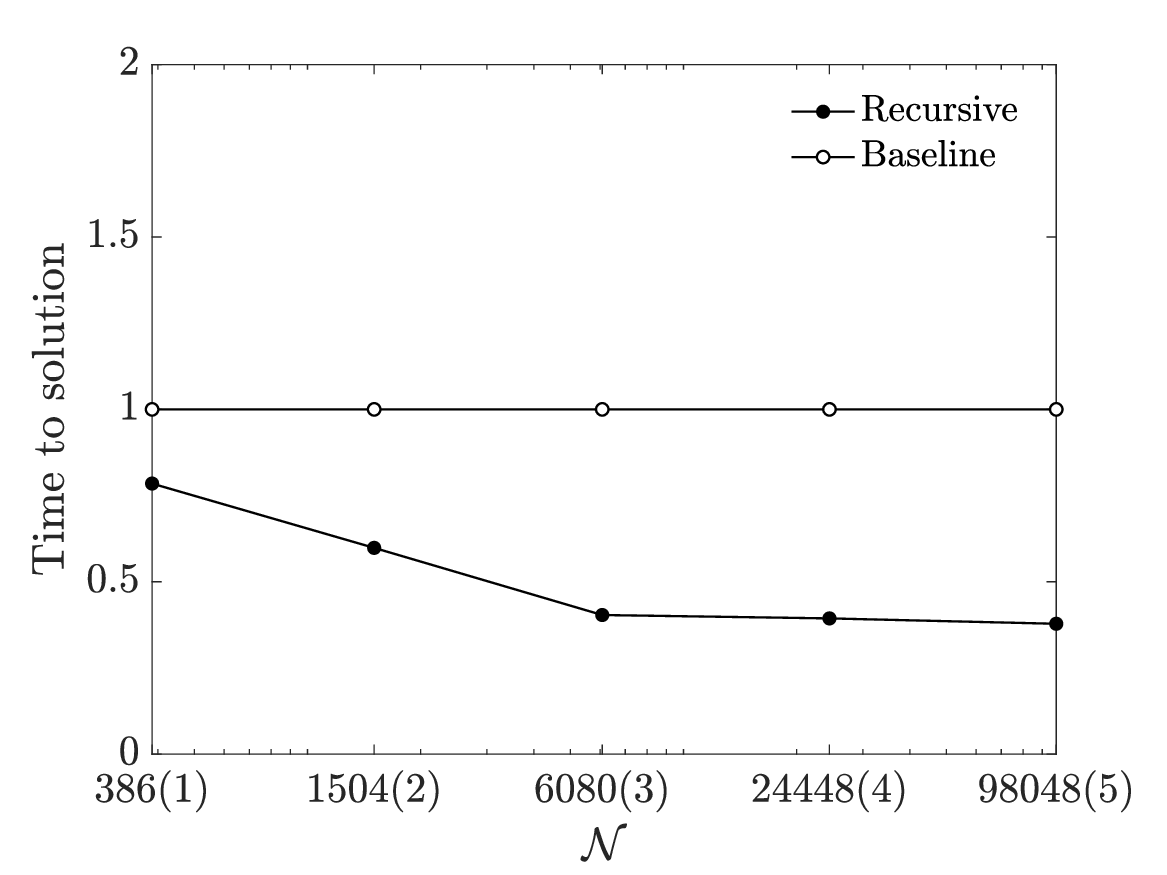}
  \caption{Linear diffusion.}
  \label{fig:ldm30}
\end{subfigure}
\hfill
\begin{subfigure}[t]{75mm}
  \includegraphics[width=75mm]{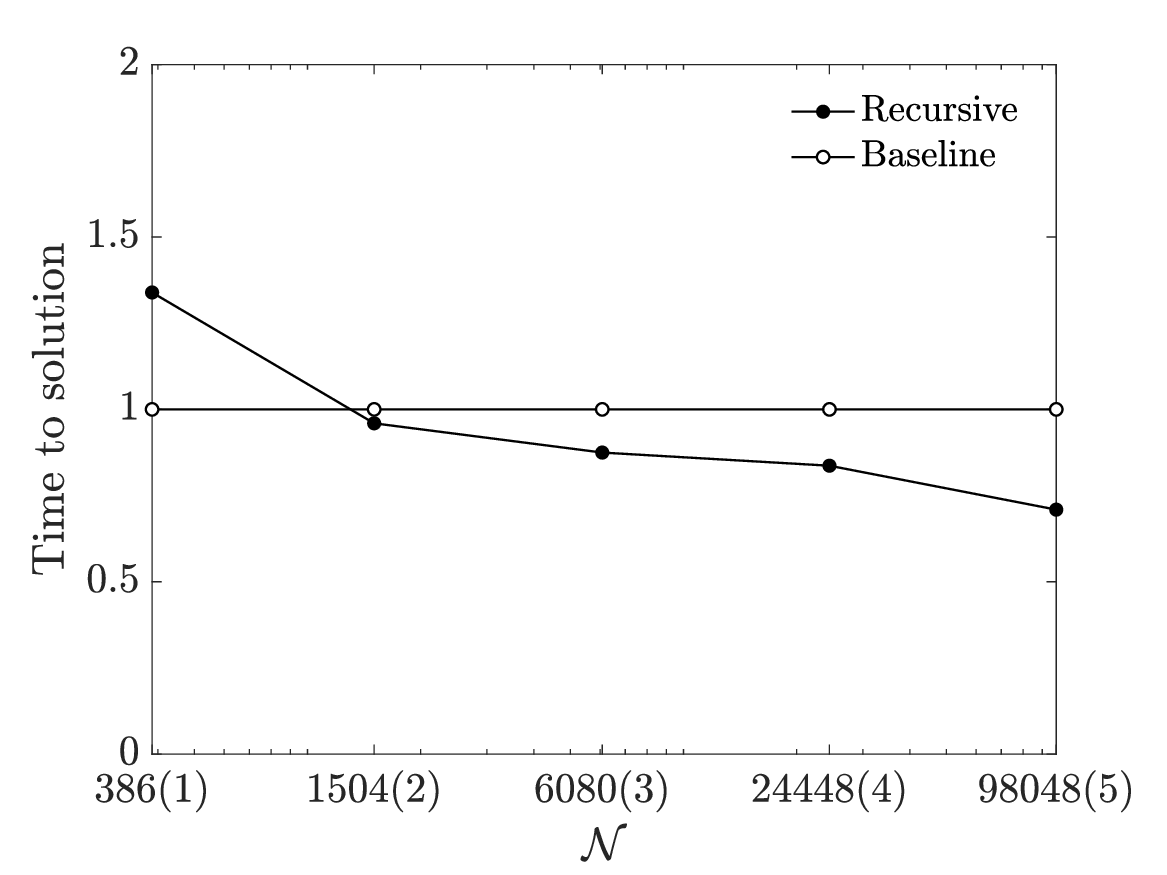}
  \caption{Convection-diffusion.}
  \label{fig:cdm30}
\end{subfigure}
\centering
\caption{Relative time comparison of the Gl-GMRES solver with and without the recursive wavelet algorithm (Algorithm \ref{alg:recAlg}) for the linear diffusion and convection-diffusion examples with $p_x=6$, $p_t=4$, $m=30(j+1)$. Level $j$ is displayed in parentheses.}
\label{fig:recursive}
\end{figure}

\begin{figure}[H]
\begin{subfigure}[t]{75mm}
  \includegraphics[width=75mm]{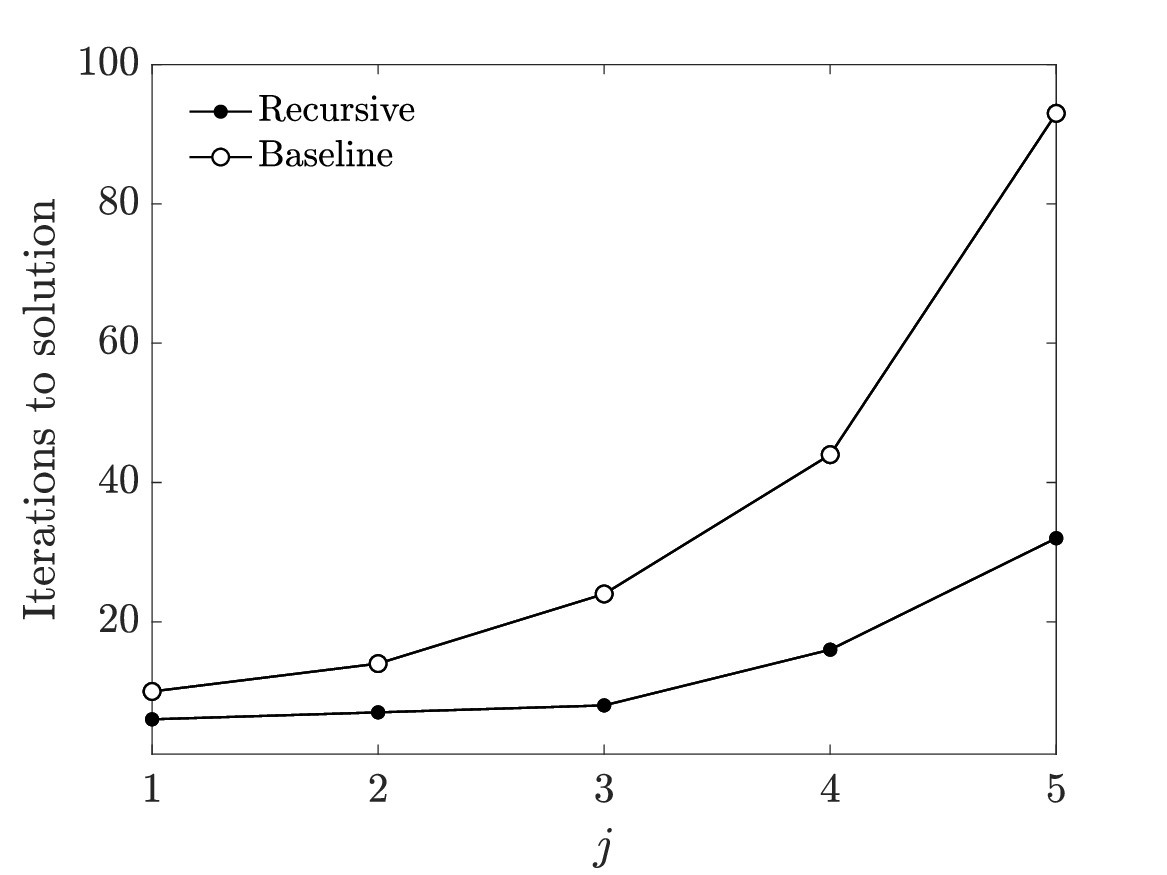}
  \caption{Linear diffusion.}
  \label{fig:ldm30It}
\end{subfigure}
\hfill
\begin{subfigure}[t]{75mm}
  \includegraphics[width=75mm]{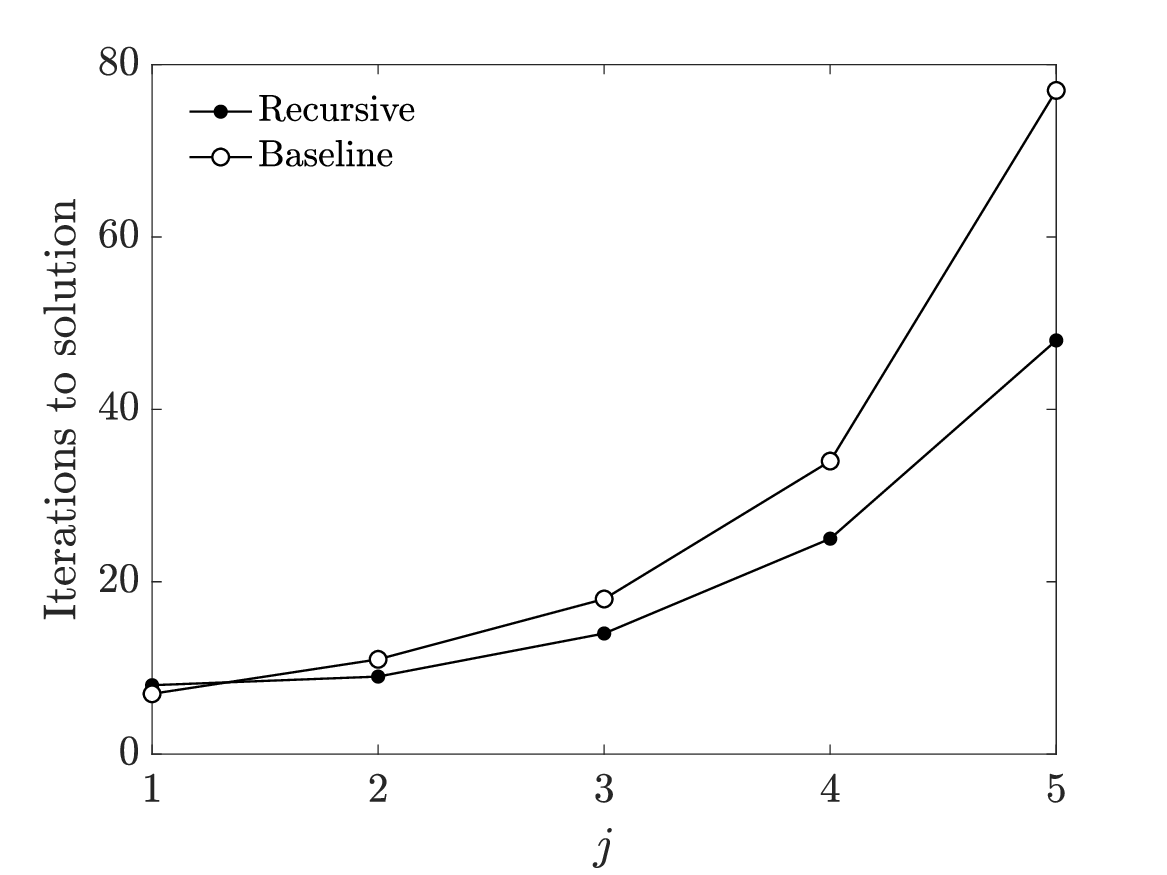}
  \caption{Convection-diffusion.}
  \label{fig:cdm30It}
\end{subfigure}
\centering
\caption{Number of Gl-GMRES iterations required to reach iterative tolerance with and without the recursive wavelet algorithm (Algorithm \ref{alg:recAlg}) for the linear diffusion and convection-diffusion examples with $p_x=6$, $p_t=4$, $m=30(j+1)$.}
\label{fig:recursiveIt}
\end{figure}

In this section, we analyze the impact of solving with Gl-GMRES using a wavelet-synthesized initial guess (Algorithm \ref{alg:recAlg}) for both linear diffusion and convection-diffusion problems. The solid horizontal line in Fig. \ref{fig:recursive} represents the ``Baseline" solution time, in which the solution is found by solving only at the desired level with an initial guess of zeros, \textit{i.e.}, there are no recursive features in the baseline  solution process. The ``Recursive" data represents the time required for the recursive technique to reach the solution relative to the baseline solution at each level. The recursive time data is cumulative for the current level and all previous levels. Therefore, the plotted recursive time data at $j = 5$  is the total time required to solve at $j = 0$, synthesize the $j=1$ initial guess, solve at $j=1$, and so on until the $j = 5$ solution is reached. As in Section \ref{sec:axbVaxxbc}, we run with restart parameter $m = 30(j+1)$, basis orders $p_x=6$ and $p_t=4$ and problem parameters $c=1$ and $\nu=0.1$. The timing data is presented as an average over multiple runs. The results for both examples are shown in Fig. \ref{fig:recursive}.
Even with the accumulated time from previous levels, the recursive solution technique is faster than the baseline method due to fewer required iterations as shown in Fig. \ref{fig:recursiveIt}. Fig. \ref{fig:recursive} shows that the number of iterations is drastically attenuated using our recursive algorithm with a better initial guess obtained from the previous level. We obtain similar results when varying the subspace size $m$.

\section{Conclusions}
\label{sec:conclusions}
In this work, we have developed a spacetime wavelet solver for linear PDEs with high-order convergence rates. We show that the spacetime wavelet formulation discretizes and solves the system in both the spatial and temporal dimensions simultaneously, obtaining an accurate solution to the resulting Sylvester matrix equation with \textit{a priori} error estimates. We enforce boundary and initial conditions and maintain well-posedness through the use of semi-orthogonal permutation matrices. The Gl-GMRES method with the embedded Modified Global Arnoldi process provides reduced memory requirements and computational effort by controlling the size of the Krylov space. The wavelet theory predicts high-order convergence for both solution and derivative approximations, which was verified with numerical experiments. By solving the systems in the Sylvester form, we are able to avoid the costly conversion to the standard $K\vec{x}=\vec{r}$ matrix form. We implement a novel wavelet-based recursive technique to speed up convergence and improve performance. Future applications of the spacetime wavelet method with Gl-GMRES would achieved improved performance with the use of carefully selected preconditioning.

\section*{Acknowledgments}
This work was supported by the Network, Cyber, and Computational Sciences Branch of the Army Research Office (ARO) under Award Number W911NF-24-1-0039. Dr. Radhakrishnan Balu served as program monitor. We would like to also thank Dr. Cale Harnish and Dr. Luke Dalessandro for their assistance through discussions on wavelet theory and computational implementation.










\newpage

\bibliography{newBib}

\begin{appendices}

\section{Matrices for Boundary/Initial Condition Enforcement} \label{sec:bcMats}
\begin{align*}
\begin{split}
    P_x &= 
    \begin{bmatrix}
        \leftarrow & \cdots & \vec{0} & \cdots & \rightarrow \\
         & & I_{n-2} & &\\
        \leftarrow & \cdots & \vec{0} & \cdots &\rightarrow
    \end{bmatrix} \in \mathbb{R}^{n \times (n-2)}\\ \\
        P_t &= 
    \begin{bmatrix}
        \leftarrow & \cdots & \vec{0} & \cdots & \rightarrow \\
         & & I_{s-1} & & \\ \\
    \end{bmatrix} \in \mathbb{R}^{s \times (s-1)} \\ \\ 
    X_D &= 
    \begin{bmatrix}
        f_{ex}(x_L,t_0) & \leftarrow & f_{ex}(x_L,\vec{t}) & \rightarrow & f_{ex}(x_L,T)  \\
        \uparrow & \ddots & \vdots & & \\
        f_{ex}(\vec{x},t_0) & \cdots & \textbf{0} & \cdots & \\
        \downarrow & & \vdots & \ddots & \\
        f_{ex}(x_R,t_0) & \leftarrow & f_{ex}(x_R,\vec{t}) & \rightarrow & f_{ex}(x_R,T)
    \end{bmatrix} \in \mathbb{R}^{n\times s}, \\ \\
    \hat{X} &= 
    \begin{bmatrix}
        f(x_1,t_1) & f(x_1,t_2) & \cdots & \cdots & f(x_1,t_{s-1})  \\
        f(x_2,t_1)  &  f(x_2,t_2) & & & \\
         \vdots & & \ddots & & \vdots \\
        f(x_{n-3},t_1) & & & \ddots &  \\
        f(x_{n-2},t_1) & f(x_{n-2}, t_2) & & f(x_{n-2},t_{s-2}) & f(x_{n-2},t_{s-1})
    \end{bmatrix} \in \mathbb{R}^{(n-2)\times (s-1)}.
    \end{split}
\end{align*}

\end{appendices}
\end{document}